\documentclass[11pt]{article}

\usepackage{amssymb,amsthm,amsmath,hyperref}

\numberwithin{equation}{section}

\newtheorem{thm}{Theorem}[section]
\newtheorem{lem}{Lemma}[section]

\newtheorem{prop}{Proposition}[section]
\theoremstyle{definition}
\newtheorem{defn}{Definition}[section]
\theoremstyle{remark}
\newtheorem{rem}{Remark}[section]

\allowdisplaybreaks

\begin{document}
\title{Global Behavior of  Spherically Symmetric
Navier-Stokes Equations with Degenerate Viscosity
    Coefficients\thanks {This
work is supported by NSFC 10571158, Zhejiang Provincial NSF of
China (Y605076) and China Postdoctoral Science Foundation
20060400335} }
\author{Mingjun Wei\thanks{E-mail: m.j.wei@126.com},
Ting Zhang\thanks{E-mail:  zhangting79@hotmail.com, Phone:
+86-571-87951860 ext 2208},
Daoyuan Fang\thanks{E-mail: dyf@zju.edu.cn, Phone: +86-571-87951860 ext 8203}\\
Department of Mathematics, Zhejiang University,\\ Hangzhou 310027,
PR China }
\date{}
\maketitle
\begin{abstract}
In this paper, we study a free boundary problem for compressible
spherically symmetric Navier-Stokes equations with a gravitational
force and degenerate viscosity coefficients. Under certain
assumptions that imposed on the initial data, we obtain the global
existence and uniqueness of the weak solution and give some
uniform bounds (with respect to time) of the solution. Moreover,
we obtain some stabilization rate estimates in $L^\infty$-norm and
weighted $H^1$-norm of the solution. The results show that such
system is stable under the small perturbations, and
could  be applied to the astrophysics.
\end{abstract}
\section{Introduction.}
We consider the compressible Navier-Stokes equations with
density-dependent viscosity in $\mathbb{R}^n(n\geq2)$, which can
be written in Eulerian coordinates as
    \begin{equation}
      \left \{
\begin{array}{l}
\partial_{\tau}\rho+\nabla\cdot(\rho \vec{u})=0, \\
\partial_{\tau}(\rho \vec{u})+\nabla\cdot(\rho \vec{u}\otimes \vec{u})+\nabla
P=\textrm{div}(\mu(\nabla \vec{u}+\nabla
\vec{u}^{\top}))+\nabla(\lambda \textrm{div}\vec{u})-\rho \vec{f},
\end{array}
      \right.\label{ss-E1.1}
    \end{equation}
in a domain $\{(\vec{\xi},\tau)\big| \vec{\xi}\in
\Omega_\tau\subset \mathbb{R}^n,\tau>0\}$, the initial conditions
and boundary conditions are
    \begin{equation}
    (\rho,\vec{u})(\vec{\xi},0)=(\rho_0,u_0)(\vec{\xi}),\ \vec{\xi}\in\Omega_0=
    \{\vec{\xi}\in\mathbb{R}^n\big| a<|\vec{\xi}|<b\},
    \end{equation}
        \begin{equation}
          \vec{u}|_{|\vec{\xi}|=a}=0,
               \  \rho|_{\vec{\xi}\in\partial \Omega_\tau\backslash \{|\vec{\xi}|=a\}}=0,\label{ss-E1.4}
                \end{equation}
where $\Omega_\tau=\psi(\Omega_0,\tau)$  and $\psi$ is the flow of
$\vec{u}$:
    \begin{equation}
      \left\{
      \begin{array}{ll}
      \partial_\tau \psi(\vec{\xi},\tau)=\vec{u}(\psi(\vec{\xi},\tau),\tau), &\vec{\xi}\in
      \mathbb{R}^n,\\
      \psi(\vec{\xi},0)=\vec{\xi}.
      \end{array}
      \right.
    \end{equation}
    Here $\rho$, $P$,
$\vec{u}=(\mathrm{u}_1,\ldots,\mathrm{u}_n)$ and $\vec{f}$ are the
density, pressure,  velocity and external force, respectively;
$\lambda=\lambda(\rho)$ and $\mu=\mu(\rho)$ are the viscosity
coefficients.

For the initial-boundary value
 problem (\ref{ss-E1.1})-(\ref{ss-E1.4}) with the spherically
symmetric initial data and external force
    $$
    (\rho,\vec{u})(\vec{\xi},0)=(\rho_0(r),u_0(r)\frac{\vec{\xi}}{r}
    ),\ \vec{ \xi}\in\Omega_0,
    $$
    $$\vec{f}=f(m,r, \tau)\frac{\vec{\xi}}{r},
    \ m(\rho,r)=\int_a^r\rho(s,\tau)s^{n-1}ds,
    \ \vec{\xi}\in\Omega_\tau,
    $$
where $r=|\vec{\xi}|=\sqrt{\xi^2_1+\cdots+\xi^2_n}$,    we are
looking for spherically symmetric solutions $(\rho,\vec{u})$:
    $$
      \rho(\vec{\xi},\tau)=\rho(r,\tau),\ \
    \vec{u}(\vec{\xi},\tau)=u(r,\tau)\frac{\vec{\xi}}{r},
    \ \vec{\xi}\in\Omega_\tau,
    $$
where  $ \Omega_\tau=\{\vec{\xi}\in\mathbb{R}^n\ \big|\
a<|\vec{\xi}|<b(\tau),b(0)=b,b'(\tau)=u(b(\tau),\tau)\} $.   Then
$(\rho,u)(r,\tau)$ is determined by
   \begin{equation}
      \left \{
\begin{array}{lll}
\partial_{t}\rho+\partial_{r}(\rho u)+\frac{n-1}{r}\rho u&=&0, \\
\rho(\partial_{t}u+u\partial_{r}u)+\partial_{r}P
&=&\partial_r[(\lambda+2\mu)(\partial_{r}u+\frac{n-1}{r}u)]-2(n-1)\partial_{r}\mu\frac{u}{r}-\rho
f,
\end{array} \right.\label{ss-E2.1}
    \end{equation}
where $(r,\tau)\in(a,b(\tau))\times(0,\infty)$, with the initial
data
    \begin{equation}
      (\rho,u)|_{\tau=0}=(\rho_0,u_0)(r),
    \  a\leq r\leq b,
    \end{equation}
 the boundary conditions
            \begin{equation}
            u|_{r=a}=0,\                  \rho|_{r=b(\tau)}
                  =0,
                  \label{ss-E2.4}
                \end{equation}
where  $b(0)=b$, $b'(\tau)=u(b(\tau),\tau)$, $\tau>0$.

 To simplify the presentation, we only consider the famous
polytropic model, i.e. $P(\rho)=A\rho^\gamma$ with $\gamma>1$ and
$A>0$ being constants.  And we assume that the viscosity
coefficients $\mu$ and $\lambda$ are proportional to
$\rho^\theta$, i.e. $\mu(\rho)=c_1\rho^\theta$ and
$\lambda(\rho)=c_2\rho^\theta$ where $c_1,c_2$ and $\theta$ are
three constants.

Additionally,  we assume the external force $ f(m, r, \tau)$
 satisfies
    \begin{equation}
    f(m,r,\tau)=f_{\infty}(m,r)+\Delta f(m,r,\tau),\label{ss-f1}
    \end{equation}
for all $m\geq0$, $r\geq a$ and     $\tau\geq0$, with
    \begin{equation}
       f_\infty(m,r)=G\frac{M_0+ m}{r^{n-1}},
       \ \Delta f(m,r,\tau)\in
    C^{1}(\mathbb{R}_+\times\mathbb{R}_+\times\mathbb{R}_+)\label{ss-f2}
        \end{equation}
         \begin{equation}
       \|\Delta f(\cdot,\cdot,\tau)\|_{L^\infty(\mathbb{R}_+\times\mathbb{R}_+)}\leq f_1(\tau),
    \          \|(\partial_r\Delta f,\partial_\tau\Delta
          f)(\cdot,\cdot,\tau)\|_{L^\infty(\mathbb{R}_+\times\mathbb{R}_+)}\leq f_2(\tau),
        \end{equation}
        \begin{equation}
      f_1\in L^\infty\cap
      L^1(\mathbb{R}_+),\ \ f_2\in {L^2}(\mathbb{R}_+),\label{ss-f4}
        \end{equation}
where $\mathbb{R}_+=[0,\infty)$,  $G>0$ is a  constant, $M_0\geq
0$ is the total mass of the solid core surrounded by the gas, and
the perturbation $\Delta f$ tends to $0$ as $\tau\rightarrow
\infty$ in some weak sense. If $M_0=0$, we ignore the
gravitational effect of the solid core. $\Delta f$ expresses the
influence of the outside gravitational force, $f_\infty$ is the
precise expression for its own gravitational force and the
gravitational force of the solid core, in the astrophysical case
(with spherical symmetry).  We study the stabilization problem of
such system, which could  be applied to the astrophysics.

Now, we consider the stationary problem,  namely
    \begin{equation}
    (P(\rho_\infty))_r=-\rho_\infty f_\infty(m(\rho_\infty,r),r)
    \label{ss-E-sym1}
    \end{equation}
in an interval $r\in(a,l_\infty)$ with the end $l_\infty$
satisfying
    \begin{equation}
     \rho_\infty(l_\infty)=0,\           \int_a^{l_\infty}\rho_\infty r^{n-1}
          dr=M:=\int_a^{b}\rho_0 r^{n-1}dr.
          \label{ss-E-sym3}
        \end{equation}
The unknown quantities are the stationary density $\rho_\infty\geq
0$ and free boundary $l_\infty>0$. If $
     \gamma>\frac{2n-2}{n}$,
 from Proposition
      \ref{ss-self-stat-ex-prop5}, we know that there exists a unique
solution $(\rho_\infty,l_\infty)$ to the stationary system
(\ref{ss-E-sym1})-(\ref{ss-E-sym3}), satisfying
$\rho_\infty(r)\sim (l_\infty^n-r^n)^{\frac{1}{\gamma-1}}$,
$(\rho_\infty)_r(r)<0$, $a<r<l_\infty$ with $l_\infty<+\infty$.

To handle the free boundary problem
(\ref{ss-E2.1})-(\ref{ss-E2.4}), it is convenient to reduce the
problem in Eulerian coordinates $(r,\tau)$ to the problem in
Lagrangian coordinates $(x,t)$ moving with the fluids, via the
transformation:
    \begin{equation}
      x=\int^r_ay^{n-1}\rho(y,\tau)dy,\ t=\tau.\label{ss-self-E1.13}
    \end{equation}
Then the fixed boundary $r=a$ and the free boundary $r=b(\tau)$
become
        $$
          x=0\ \textrm{ and }x=\int^{b(\tau)}_ay^{n-1}\rho(y,\tau)dy
          =\int^{b}_ay^{n-1}\rho_0(y)dy=M,
        $$
where $M$ is the total mass initially. So that the region
$\{(r,\tau)|a\leq r\leq b(\tau),\tau\geq0\}$ under consideration
is transformed into the region $\{(x,t)|0\leq x\leq M, t\geq0\}$,
and the function $m(\rho,r)$ becomes $x$.
    Under the coordinate transformation (\ref{ss-self-E1.13}), the equations
    (\ref{ss-E2.1})-(\ref{ss-E2.4}) are transformed into
   \begin{equation}
      \left\{\begin{array}{lll}
      \partial_t \rho(x,t)&=&-\rho^2\partial_x(r^{n-1} u),\\
      \\
           \partial_t u(x,t)&=&r^{n-1}\left\{\partial_x[\rho(\lambda+2\mu)\partial_{x}(r^{n-1}u)-P]
           -2(n-1)\frac{u}{r}\partial_x\mu\right\}-f(x,r,t),\\ \\
      r^n(x,t)&=&a^n+n\int^x_0\rho^{-1}(y,t)dy,
     \end{array}
      \right.\label{ss-E}
    \end{equation}
 where $(x,t)\in(0,M)\times(0,\infty)$, with the initial data
    \begin{equation}
      (\rho,u)|_{t=0}=(\rho_0,u_0)(x),
        r|_{t=0}=r_0(x)=\left(a^n+
        n\int^x_0\rho^{-1}_0(y)dy
        \right)^\frac{1}{n},
    \end{equation}
and  the boundary conditions
            \begin{equation}
            u(0,t)=0,\                 \rho|_{x=M}=0,
                  \ t>0.
                  \label{ss-Efd}
                \end{equation}

 It is standard that if we can solve the problem
(\ref{ss-E})-(\ref{ss-Efd}), then the free boundary problem
(\ref{ss-E1.1})-(\ref{ss-E1.4}) has a solution.

From (\ref{ss-E-sym1})-(\ref{ss-E-sym3}), it is easy to see that
$\rho_{\infty}(x)$ is the solution to the stationary system,
    \begin{equation}
      Ar_\infty^{n-1}(\rho_\infty^\gamma)_x=-f_{\infty}(x,r_{\infty}),
      \ r_{\infty}^n(x)=a^n+n\int^x_0\rho^{-1}_{\infty}(y)dy,
      \ x\in(0,M),\label{ss-E1.17}
    \end{equation}
        \begin{equation}
          \rho_\infty(M)=0.\label{ss-E1.18}
        \end{equation}

The results in \cite{hoff,xin} show that the compressible
Navier-Stokes system with the constant viscosity coefficient have
the singularity at the vacuum.  Considering the modified
Navier-Stokes system in which the viscosity coefficient depends on
the density, Liu-Xin-Yang\cite{liu} proved that such system is
local well-posedness. It is motivated by the physical
consideration that in the derivation of the Navier-Stokes
equations from the Boltzmann equation through the Chapman-Enskog
expansion to the second order, cf. \cite{Grad}, the viscosity is a
function of the temperature. If we consider the case of isentropic
fluids, this dependence is reduced to the dependence on the
density function.

 Since $n\geq2$ and the viscosity coefficient $\mu$ depends on $\rho$,
 the  nonlinear term  $2(n-1)\frac{1}{r}u\partial_x\mu$ in (\ref{ss-E})$_2$ makes the analysis
  significantly different from  the one-dimensional
case \cite{liu,okada1,vong,yang3,fang06}.  When
$\mu\geq\underline{\mu}>0$ and $\rho_0\geq \underline{\rho}>0$,
authors in \cite{Ducomet2005,Zlotnik2005} obtained the existence,
uniqueness and global behavior of the solution for  compressible
spherically symmetric Navier-Stokes equations with a external
pressure and without the nonlinear term
$2(n-1)\frac{1}{r}u\partial_x\mu$. Following the ideas in
\cite{Zlotnik2005}, we can obtain the existence and uniqueness
results for the stationary problem in Section
\ref{ss-self-sec1.5}. In this paper, since viscosity coefficients
and density will degenerate at the free boundary, we need to use
the weighted function $(M-x)$ to control the lower bound of the
density in Section \ref{ss-sec2}.

Considering the system (\ref{ss-E})-(\ref{ss-Efd}) with a general
external force, Chen-Zhang obtained the local existence and
uniqueness of the solution in \cite{ChenP}. In this paper, when
the initial data $(\rho_0,u_0,r_0)$ are close to the stationary
state $(\rho_\infty,0,r_\infty)$, we will obtain some appropriate
\textit{a priori} estimates and prove that the maximum existence
time $T^*=\infty$. The difficulty of this problem is to obtain the
lower bound of the density $\rho$. The key ideas are using the
classical continuity method and the result of \textbf{Claim 1}. In
\textbf{Claim 1}, we want to prove that there is a small positive
constant $\epsilon_1$, such that, for any $T>0$, if
     \begin{equation}
      I(t)=\|g(\cdot,t)-g_\infty\|_{L^\infty}     \leq
      2\epsilon_1,\ \forall\ t\in[0,T],
    \end{equation}
where $g(x,t)=(M-x)^{-\frac{1}{\gamma}}\rho(x,t)$ and
$g_\infty(x)=(M-x)^{-\frac{1}{\gamma}}\rho_\infty(x)$, then
    \begin{equation}
    I(t)\leq \epsilon_1, \ \forall\ t\in[0,T].
    \end{equation}
Using the energy method and induction method, we can estimate the
weighted $L^2-$norm of $g-g_\infty$ in Lemma \ref{ss-L3.8}. In
such process (see Lemmas \ref{ss-L3.5-0}-\ref{ss-L3.7}), we use
the weight function $(1+t)^{\alpha}$ (with $\alpha=-\frac{5}{8}$)
to remedy the disadvantage of the nonlinear term
$2(n-1)\frac{1}{r}u\partial_x\mu$, and use the induction method to
increase $\alpha$ to $-\epsilon_2$. Then, by the reduction to
absurdity, we can finish the proof of \textbf{Claim 1} in Lemma
\ref{ss-Lem3.4}.

 Our results
show that: such system is stable under the small perturbations,
 does not develop vacuum states or concentration states
for all time, and the interface
 $\partial\Omega_\tau$ propagates with
finite speed.

 The assumptions on $c_1$, $c_2$,
$\theta$, $\gamma$ and initial data can be stated as follows:
    \begin{description}
\item [(A1)]  $ \gamma>\frac{2n-2}{n}$, $\theta\in
(0,\gamma-1)\cap (0,\frac{\gamma}{2}]$, $c_1$ and $c_2$ satisfy
that
        $
    c_1>0$ and $2c_1+nc_2>0;
        $
\item [(A2)] $N_1(M-x)^{1/\gamma}\leq\rho_0\leq
N_2(M-x)^{1/\gamma}$, with some positive constants $0<N_1<
 N_2$, and
$(M-x)^{1-\frac{\theta}{\gamma}}(\rho^\theta_0)_x^2\in
L^1([0,M])$, $ (\rho_0^{\gamma})_x\in L^{2}([0, M])$; \item
[(A3)]$u_0\in L^2([0,M])$, $\rho_0^{\frac{\theta+1}{2}}(u_0)_x\in
L^2([0,M])$, $u_0(0)=0$,
    \begin{equation}
    \left((2c_1+c_2)\rho_{0}^{\theta+1}(r_{0}^{n-1}u_{0})_{x}\right)_{x}
    -2c_1(n-1)\frac{u_0}{r_0}(\rho_0^{\theta})_x\in
    L^{2}([0, M]).\label{ss-1.16}
    \end{equation}
\end{description}

Under the above assumptions (A1)-(A3), we will prove the existence
of the global weak solution to the initial-boundary value problem
(\ref{ss-E})-(\ref{ss-Efd}) in the sense of the following
definition.
\begin{defn} A pair of functions $(\rho,u,r)(x,t)$ is called a
global weak solution to the initial boundary value problem
(\ref{ss-E})-(\ref{ss-Efd}) if for any $T>0$,
    \begin{equation*}
    \rho,u \in L^\infty([0,M]\times[0,T])\cap C^1([0,T];L^2([0,M])),
    \end{equation*}
        \begin{equation*}
          r\in C^{1}([0,T];L^\infty([0,M])),
        \end{equation*}
        \begin{equation*}
     \rho^{-1},\  (r^{n-2}u)_x,\   (r^{n-1})_x\in L^\infty([0,T];L^{1}([0,M])),
        \end{equation*}
 and       \begin{equation*}
         \rho^{1+\theta}(r^{n-1}u)_x\in L^\infty([0,M]\times[0,T])\cap C^\frac{1}{2}([0,T];L^2([0,M])).
        \end{equation*}
Furthermore, the following equations hold:
    \begin{equation*}
    \rho_t+\rho^2(r^{n-1}u)_x=0,\
    \rho(x,0)=\rho_0(x),\ a.e.
    \end{equation*}
    \begin{equation*}
      r_t=u,\ r(x,0)=r_0(x),\ r^{n}(x, t)=a^{n}+n\int_{0}^{x}\rho^{-1}(y, t)dy,\ a.e.
    \end{equation*}
        \begin{eqnarray*}
    &&\int^\infty_0\int^M_0[u\psi_t+(P-\rho(\lambda+2\mu)(r^{n-1}
    u)_x)(r^{n-1}\psi)_x\nonumber\\
    &&+2(n-1)\mu(r^{n-2}u\psi)_x-f(x,r,t)\psi]dxdt
        +\int^M_0u_0(x)\psi(x,0)dx=0,
        \end{eqnarray*}
for any test function $\psi(x,t)\in C^\infty_0(\Omega)$ with
$\Omega=\{(x,t)\ \big|\ 0< x\leq M,\ t\geq 0\}$, and
$$
\lim_{\epsilon\rightarrow0}\frac{1}{\epsilon}\int^\epsilon_0
udx=\lim_{\epsilon\rightarrow0}\frac{1}{\epsilon}\int^1_{1-\epsilon}
\rho dx=0.
$$
\end{defn}

In what follows, we always use $C$($C_i$) to denote a generic
positive constant depending only on $\gamma$, $\theta$, $f_1$,
$f_2$ and the initial data, independent of the given time T.

 We
now state the main theorems in this paper.
\begin{thm}\label{ss-thm} Under the conditions (\ref{ss-f1})-(\ref{ss-f4}) and [A1-A3],  there exists a
constant $\epsilon_0>0$, such that if
    \begin{equation}
             \|f_1\|^2_{L^\infty}+\int^\infty_0(1+t)f_1^2(t)dt\leq \epsilon_0^2,\label{ss-epsilon1}
    \end{equation}
    \begin{equation}
      \|u_0\|_{L^2}+\|(M-x)^{-\frac{1}{\gamma}}(\rho_0-\rho_\infty)\|_{L^\infty}\leq
      \epsilon_0,\label{ss-epsilon2}
    \end{equation}
 then the system
(\ref{ss-E})-(\ref{ss-Efd}) has a unique global weak solution
$(\rho,u,r)$ satisfying
\begin{equation}
          C^{-1}(M-x)^\frac{1}{\gamma}\leq
          \rho(x,t)\leq           C(M-x)^\frac{1}{\gamma},
         \label{ss-E1.32}
        \end{equation}
  \begin{equation}
    r(x,t)\in [a,C],\label{ss-E1.32-7}
  \end{equation}
               \begin{equation}
                  \int^M_0(M-x)^{1-\frac{\theta}{\gamma}}
                  (\rho^\theta-\rho_\infty^\theta)_x^2dx\leq
                  C,
                \end{equation}
  and
         \begin{equation}
             \left\|u(\cdot,t)
             \right\|_{L^\infty}+ \|[\rho(r^{n-1}u)_x](\cdot,t)\|_{L^\infty}\leq
              C,\label{ss-E1.46}
            \end{equation}
for all $t\geq0$ and $x\in[0,M]$. Furthermore, if
$(1+t)^{\frac{2(\gamma+\theta)}{\gamma+\theta+1}}f_1^2(t)\in
L^1(\mathbb{R}_+)$, for any $\eta>0$, we have
        \begin{equation}
      \int^M_0\left\{u^2+(M-x)^{\frac{\gamma-1}{\gamma}}(g-g_\infty)^2
      +(r-r_\infty)^2\right\}dx\leq C_{\eta}(1+t)^{-\frac{2(\gamma+\theta)}{\gamma+\theta+1}+\eta},\label{ss-E1.32-9}
        \end{equation}
            \begin{equation}
              \int^M_0\left(\rho^{\theta-1}
              u^2+\rho^{\theta+1}u_x^2
              \right)dx\leq
              C_{\eta}(1+t)^{-\frac{2(\gamma+\theta)}{\gamma+\theta+1}+\eta},\label{ss-E1.44}
            \end{equation}
                   \begin{equation}
              \int^M_0(M-x)^{2-\frac{2\theta}{\gamma}}(\rho^\theta-\rho_\infty^\theta)_x^2dx\leq
              C_\eta(1+t)^{-\frac{\gamma+\theta-1}{\gamma+\theta+1}+\eta},
            \end{equation}
and
            \begin{equation}
            \|\rho^{\gamma}(\cdot,t)-\rho_\infty^{\gamma}(\cdot)\|_{L^\infty}\leq
             C_\eta(1+t)^{-\frac{3\gamma+3\theta-1}{4(\gamma+\theta+1)}+\frac{\eta}{2}},
              \  \|\rho^\frac{\gamma+\theta}{2}(\cdot,t)-\rho_\infty^\frac{\gamma+\theta}{2}(\cdot)\|_{L^\infty}\leq
             C_{\eta}(1+t)^{-\frac{\gamma+\theta}{2(\gamma+ \theta+1)}+\frac{\eta}{4}},
             \label{ss-E1.43-0}
            \end{equation}
             \begin{equation}
             \|u(\cdot,t)\|_{L^\infty}\leq C_{\eta}
             (1+t)^{-\frac{\gamma+\theta}{\gamma+ \theta+1}+\frac{\eta}{2}},\label{ss-E1.43}
            \end{equation}
for all $t\geq0$, where $C_\eta$ is a positive constant depending
on $\eta$.
\end{thm}

\begin{rem}
  The uniqueness of the solution in Theorems \ref{ss-thm} or
  \ref{ss-localthm} means that: if $(\rho_1,u_1,r_1)$ and
  $(\rho_2,u_2,r_2)$ are two solutions to the
system (\ref{ss-E})-(\ref{ss-Efd}) with the same initial data
$(\rho_{0},u_{0},r_{0})$, and satisfy regularity conditions in the
theorem, then we have $(\rho_1,u_1,r_1)=(\rho_2,u_2,r_2)$.
\end{rem}

\begin{rem}
  In particular, the viscosity of the gas is proportional to the
square root of the temperature for the hard sphere model (as
pointed out in \cite{okada1,yang3}), and the relation between
$\theta$ and $\gamma$ is
    $$
      \theta=\frac{\gamma-1}{2}.
    $$
Our condition (A1) covers it. Since the Navier-Stokes system with
constant viscosity coefficients has the singularity at the vacuum
\cite{hoff,xin}, we assume $\theta>0$ in (A1).
\end{rem}
\begin{rem}
There are no smallness assumptions on
$\|(M-x)^{1-\frac{\theta}{\gamma}}(\rho^\theta_0)_x^2\|_{L^1}$ and
$\|\rho_0^\frac{1+\theta}{2}( u_0)_x\|_{L^2}$. Moreover,
$\epsilon_0$ do not depend on
$\|(M-x)^{1-\frac{\theta}{\gamma}}(\rho^\theta_0)_x^2\|_{L^1}$ and
$\|\rho_0^\frac{1+\theta}{2}( u_0)_x\|_{L^2}$.
\end{rem}
\begin{rem}
Considering the no vacuum system in an exterior domain in
$\mathbb{R}^3$, Kobayashi-Shibata\cite{Kobayashi} obtained
$\|(\rho-\rho_\infty,u)(\cdot,t)\|_{L^2}\lesssim(1+t)^{-\frac{3}{4}}$
and
$\|(\rho-\rho_\infty,u)(\cdot,t)\|_{L^\infty}\lesssim(1+t)^{-\frac{3}{2}}$
when $\rho_\infty$ is a positive constant. Considering the no
vacuum system in $\mathbb{R}^n$ ($n\geq3$),
Ukai-Yang-Zhao\cite{Ukai} obtained
$\|(\rho-\rho_\infty,u)(\cdot,t)\|_{L^2\cap L^\infty}\leq
C(1+t)^{-\frac{n}{4}+\epsilon}$ when $\rho_\infty$ is  close to a
positive constant. Considering the one dimensional system with a
degenerate viscosity coefficient, we\cite{fang06} obtained
$\|(\rho^\gamma-\rho_\infty^{\gamma}, u)(\cdot,t)\|_{L^\infty}\leq
C(1+t)^{-\frac{1}{2}}$ when the external force $f_\infty$ is close
to a positive constant $N_0$ and the stationary density
$\rho_\infty$ is close to
$\left(\frac{N_0(M-x)}{A}\right)^\frac{1}{\gamma}$. Since
$\frac{\gamma+\theta}{\gamma+\theta+1}>\frac{1}{2}$ and
$\frac{3\gamma+3\theta-1}{4(\gamma+\theta+1)}>\frac{1}{2}$ (if
$\gamma+\theta>3$), it is easy to see that our results in Theorem
\ref{ss-thm} are better than the results in \cite{fang06}.  Using
 similar arguments as that in Theorem \ref{ss-thm}, we can also
   obtain  similar results in one-dimensional case which are better than the results in \cite{fang06}. For example,  the stabilization rate
   estimate $\| (\rho^\gamma-\rho_\infty^{\gamma},
u)(\cdot,t)\|_{L^\infty}\leq C(1+t)^{-\frac{1}{2}}$ in
\cite{fang06} can be replaced by
    $\| (\rho^\gamma-\rho_\infty^{\gamma},
u)(\cdot,t)\|_{L^\infty}\leq C_{\eta}(1+t)^{-
              \frac{\gamma+\theta}{\gamma+ \theta+1}+\frac{\eta}{2}}  $,
             $\|\rho^\frac{\gamma+\theta}{2}(\cdot,t)-\rho_\infty^\frac{\gamma+\theta}{2}(\cdot)\|_{L^\infty}\leq
             C_{\eta}(1+t)^{-\frac{\gamma+\theta}{2(\gamma+ \theta+1)}+\frac{\eta}{4}}$,
              $t\geq0$, for any $\eta>0$.
\end{rem}
\begin{rem}
Since that the information of the dimension $n$ mainly appear in
the index of the radii $r$, and we only consider the system with a
solid core $r\geq a>0$ in this paper, our results can not show the
effect of the dimension. In \cite{fang06-2}, we studied the global
behavior of the solution to the similar problem with a positive
external pressure and without a solid core, and obtained the
stabilization rate estimates for the solution of exponential type.
The admissible range of the parameter $\frac{\lambda}{\mu}$
depends on the dimension $n$ in \cite{fang06-2}. We will study the
system without a solid core $r\geq0$ and with degenerate viscosity
coefficients in the future, and guess that stabilization rate
estimates of the solution  can not better than the results in
Theorem \ref{ss-thm}.
\end{rem}
\begin{thm}\label{ss-thm-2} (Continuous Dependence)
For each $i=1,2$, let $(\rho_i,u_i,r_i)$ be the solution to the
system (\ref{ss-E})-(\ref{ss-Efd}) with the initial data
$(\rho_{0i},u_{0i},r_{0i})$, which satisfies regularity conditions
in Theorem \ref{ss-thm}. Then, we have
   \begin{eqnarray*}
&&\int^M_0\left((u_1-u_2)^2+\rho_1^{1-\theta}\rho_2^{2\theta-4}(\rho_1-\rho_2)^2+\rho_1^{\theta}\rho_2^{-1}(r_1-r_2)^2\right)dx\\
&\leq&Ce^{Ct}\int^M_0\left((u_{01}-u_{02})^2+\rho_{01}^{1-\theta}\rho_{02}^{2\theta-4}(\rho_{01}-\rho_{02})^2+\rho_{01}^{\theta}\rho_{02}^{-1}(r_{01}-r_{02})^2\right)dx,
\end{eqnarray*}
for all $t\geq0$.
\end{thm}

\begin{rem}
Using  similar arguments as that in  \cite{ChenP}, we can easily
obtain such continuous dependence of the solution on the initial
data, and omit the detail.
\end{rem}

\begin{rem}
If we ignore the influence of self-gravitation, i.e.
  assume $f_\infty(m,r)=G\frac{M_0}{r^{n-1}}$ with $M_0>0$, then  we can also
  obtain the same results in Theorem \ref{ss-thm}-\ref{ss-thm-2}.
\end{rem}

We now briefly review the previous works in this direction.
 For the related free boundary problem of one-dimensional isentropic
fluids with density-dependent viscosity (like
$\mu(\rho)=c\rho^\theta$), see \cite{liu,okada1,
vong,yang3,fang06} and the references therein.  For the related
stabilization rate estimates of 1-D free boundary problem, see
\cite{Ducomet2005-2,Matsumura,Straskraba,fang06} etc.. For the
spherically symmetric solutions of the Navier-Stokes equations
with
 a free boundary, see \cite{Chen2002,Ducomet2005,
MatusuNecasova,okada93,fang06-2,Zlotnik2005} etc..  Also see
Bresch-Desjardins\cite{Bresch2006}, Lions\cite{Lions} and
Vaigant-Kazhikhov\cite{Vaigant} for multidimensional isentropic
fluids.

The rest of this paper is organized as follows. First, we obtain
the existence and uniqueness of the solution to the stationary
problem in Section \ref{ss-self-sec1.5}. In Section \ref{ss-sec2},
we will prove some \textit{a priori} estimates, and extend the
local solution in \cite{ChenP} to the global solution in time. In
Section \ref{ss-Sec4}, we obtain the stabilization rate estimates
of the solution.

\section{The stationary problem}\label{ss-self-sec1.5}
Zlotnik-Ducomet\cite{Zlotnik2005} obtained the existence of the
positive solution to the stationary problem with a positive
external pressure. Using  similar arguments as that in
\cite{Zlotnik2005}, we can obtain the following results for the
stationary problem without a external pressure. We start with a
proof of the existence of a non-negative solution to the
Lagrangian stationary problem.

\begin{prop}\label{ss-self-stat-ex-prop1}
If
    \begin{equation}
    \gamma>\frac{2n-2}{n}\label{ss-self-E2.1}
    \end{equation}
or
    \begin{equation}
      \gamma=\frac{2n-2}{n}
      \ \textrm{ and }\ \left(\frac{n\gamma}{(\gamma-1)}M^\frac{\gamma-1}{\gamma}\right)^{\frac{2n-2}{n}}
    < \frac{G}{A}\left(\frac{M}{2}+M_0\right),\label{ss-self-E2.2}
    \end{equation}
or
    \begin{equation}
 n>2,\    1<\gamma<\frac{2n-2}{n}
      \ \textrm{ and }\
              \delta_6^\gamma\left(a^n+\frac{n\gamma}{\delta_6(\gamma-1)}M^\frac{\gamma-1}{\gamma}
              \right)^{\frac{2n-2}{n}}
    \leq \frac{G}{A}\left(\frac{M}{2}+M_0\right),
    \end{equation}
 where $\delta_6=a^n(1-\frac{\gamma
 n}{2n-2})\frac{2n-2}{\gamma-1}M^\frac{\gamma-1}{\gamma}$,
 then the Lagrangian stationary problem
(\ref{ss-E1.17})-(\ref{ss-E1.18}) has a non-negative solution
$\rho_\infty\in W^{1,\beta}([0,M])$ satisfying
$C^{-1}(M-x)^\frac{1}{\gamma}\leq \rho_\infty(x)\leq
C(M-x)^\frac{1}{\gamma}$, where
$\beta\in[1,\frac{\gamma}{\gamma-1})$ is a constant.
\end{prop}
\begin{proof}
  We introduce the nonlinear operator
    $$
    I:K\rightarrow W^{1,\beta}([0,M]),
    $$
where
    $K=\left\{f\in
    C([0,M])\left| f\geq0, \|\frac{(M-x)^{\frac{1}{\gamma}}}{f(x)}\|_{L^\infty}<\infty,
    \ \|\frac{f(x)}{(M-x)^{\frac{1}{\gamma}}}\|_{L^\infty}<\infty,
    \right.
    \right\}$, by setting
        $$
        I(f)(x)=\left(\frac{\int^M_xG\frac{M_0+y}{r_f^{2n-2}(y)}dy}{A}
        \right)^\frac{1}{\gamma}
    \textrm{with}\ r_f^n(x)=a^n+n\int^x_0f^{-1}(y)dy,\
    x\in[0,M].
    $$
We can restate the problem (\ref{ss-E1.17})-(\ref{ss-E1.18}) as
the fixed-point problem
    \begin{equation}
    \rho_\infty=I (\rho_\infty).\label{ss-self-E2.3}
    \end{equation}

For all $f\in K_{\delta_1,\delta_2}=\left\{ f\in
K\left|\delta_1(M-x)^{\frac{1}{\gamma}}\leq f(x)\leq
\delta_2(M-x)^\frac{1}{\gamma}\right.\right\}$ with
$0<\delta_1\leq\delta_2<\infty$, we have
        $$
        a^n\leq r_\infty^n(x)\leq
        a^n+\frac{n\gamma}{\delta_1(\gamma-1)}M^\frac{\gamma-1}{\gamma}:=B^n
        $$
and
    $$
    \left(\frac{G(\frac{M}{2}+M_0)}{AB^{2n-2}}
    \right)^\frac{1}{\gamma}(M-x)^\frac{1}{\gamma}
    \leq I(f)(x)
    \leq     \left(\frac{G(M+M_0)}{Aa^{2n-2}}
    \right)^\frac{1}{\gamma}(M-x)^\frac{1}{\gamma}, x\in[0,M].
    $$

If $\gamma>\frac{2n-2}{n}$, then $I (K_{\delta_3,\delta_4})\subset
K_{\delta_3,\delta_4}$, where
$\delta_4=\left(\frac{G(M+M_0)}{Aa^{2n-2}}
    \right)^\frac{1}{\gamma}$ and $\delta_3$ is a positive constant satisfying
    $\delta_3^\gamma(a^n+\frac{n\gamma}{\delta_3(\gamma-1)}M^\frac{\gamma-1}{\gamma})^{\frac{2n-2}{n}}
    \leq \frac{G}{A}(\frac{M}{2}+M_0)$. And one can immediately verify that $I$ is a
    compact operator on $K_{\delta_3,\delta_4}$. Since $K_{\delta_3,\delta_4}$  is a convex
 closed bounded non-empty subset of $C([0,M])$, the problem
 (\ref{ss-self-E2.3}) has a solution $\rho\in K_{\delta_3,\delta_4}$ by Schauder's fixed point theorem.

Similarly, if $\gamma=\frac{2n-2}{n}$ and
$(\frac{n\gamma}{(\gamma-1)}M^\frac{\gamma-1}{\gamma})^{\frac{2n-2}{n}}
    < \frac{G}{A}(\frac{M}{2}+M_0)$,
then $I (K_{\delta_5,\delta_4})\subset K_{\delta_5,\delta_4}$,
where
$\delta_5=a^{-n}\left[\left(\frac{G}{A}(\frac{M}{2}+M_0)\right)^\frac{n}{2n-2}-\frac{n\gamma}{(\gamma-1)}M^\frac{\gamma-1}{\gamma}
\right]$, and problem
 (\ref{ss-self-E2.3}) has a solution $\rho\in K_{\delta_5,\delta_4}$.

Similarly, if $n>2$, $1<\gamma<\frac{2n-2}{n}$ and
    $         \delta_6^\gamma(a^n+\frac{n\gamma}{\delta_6(\gamma-1)}M^\frac{\gamma-1}{\gamma})^{\frac{2n-2}{n}}
    \leq \frac{G}{A}(\frac{M}{2}+M_0)$,
then $I( K_{\delta_6,\delta_4})\subset K_{\delta_6,\delta_4}$,
and problem
 (\ref{ss-self-E2.3}) has a solution $\rho\in K_{\delta_6,\delta_4}$.
\end{proof}

Similar to \cite{Zlotnik2005},
 We say a stationary solution $(\rho_\infty,r^n_\infty)$ is \textit{statically
stable} if
    \begin{eqnarray}
        J[W]&:=&
      \int^M_0\left(\gamma A\rho_\infty^{1+\gamma}W_x^2
    -(2n-2)G(M_0+x)r_\infty^{2-3n} W^2
    \right)dx\nonumber\\
      &\geq &\delta_{7} \int^M_0\left((M-x)^\frac{1+\gamma}{\gamma}W_x^2
    + W^2
    \right)dx,\label{ss-selfE2.8}
    \end{eqnarray}
for some $\delta_7>0$ and all $$W\in K_1=\left\{f\in
C([0,M])\left| f(0)=0, \left\|(M-x)^{\frac{1}{\gamma}}
f'(x)\right\|_{L^\infty}<\infty,
\left\|\frac{1}{(M-x)^{\frac{1}{\gamma}}
f'(x)}\right\|_{L^\infty}<\infty\right.\right\}.$$

Now, the static potential energy  takes the following form:
    \begin{equation}
        S[V]=\int^M_0\left(\frac{A}{\gamma-1}(V_x)^{1-\gamma}
        +\int^V_{\frac{a^n}{n}}G(M_0+x)(nh)^{\frac{2-2n}{n}}dh\label{ss-self-E2.9-1}
    \right)dx.
    \end{equation}
We call  $V\in K_2=\left\{f\in C([0,M])\left| f(0)=\frac{a^n}{n},
\left\|(M-x)^{\frac{1}{\gamma}} f'(x)\right\|_{L^\infty}<\infty,
\left\|\frac{1}{(M-x)^{\frac{1}{\gamma}}
f'(x)}\right\|_{L^\infty}<\infty\right.\right\}$ is a point of
\textit{local quadratic minimum of $S$} if
    \begin{equation}
      S[V+W]-S[V]\geq \delta_8\int^M_0\left((M-x)^\frac{1+\gamma}{\gamma}W_x^2
    + W^2
    \right)dx,\label{ss-self-E2.9}
    \end{equation}
for all $W\in K_1$ and
$\|(M-x)^{\frac{1}{\gamma}}W_x\|_{L^\infty([0,M])}+\|W\|_{L^\infty}\leq
\delta_9$, for some $\delta_8>0$ and $\delta_9>0$.

\begin{prop}\label{ss-self-stat-ex-prop3}
If $\gamma>\frac{2n-2}{n}$ and $\rho_\infty$ is a solution of the
problem (\ref{ss-E1.17})-(\ref{ss-E1.18}) satisfying
$\rho_\infty\in W^{1,\beta}([0,M])$ and
$C^{-1}(M-x)^\frac{1}{\gamma}\leq \rho_\infty(x)\leq
C(M-x)^\frac{1}{\gamma}$, then we have (\ref{ss-selfE2.8})    and
(\ref{ss-self-E2.9}) hold with $V=V_\infty=\frac{r_\infty^n}{n}$.
\end{prop}
\begin{proof}
  From $r_\infty\geq a$, $(A\rho_\infty^\gamma)_x=-G\frac{M_0+x}{r_\infty^{2n-2}}$ and $(r_\infty^n)_x=n\rho_\infty^{-1}$,
  using integration by parts,  we
  have
    \begin{eqnarray*}
    J[W]
            &=&\int^M_0\left(\gamma A\rho_\infty^{1+\gamma}W_x^2
    +(2n-2)A(\rho_\infty^\gamma)_xr_\infty^{-n} W^2
    \right)dx\nonumber\\
            &=&\int^M_0\left(\gamma A\rho_\infty^{1+\gamma}W_x^2
    -2(2n-2)A\rho_\infty^\gamma r_\infty^{-n} WW_x
    \right.\nonumber\\
    &&\left.+n(2n-2)A\rho_\infty^{\gamma-1}r_\infty^{-2n} W^2
    \right)dx,
    \ \textrm{for all }\ W\in K_1.
    \end{eqnarray*}
If $\gamma>\frac{2n-2}{n}$, we have
    $$
        J[W]\geq C^{-1} \int^M_0\left((M-x)^\frac{1+\gamma}{\gamma}W_x^2
    +(M-x)^\frac{\gamma-1}{\gamma} W^2
    \right)dx.
    $$
From $r_\infty\geq a$ and
$(A\rho_\infty^\gamma)_x=-G\frac{M_0+x}{r_\infty^{2n-2}}$, using
integrating by parts and the Cauchy-Schwarz inequality, we have
\begin{eqnarray}
      \int^M_0W^2dx&\leq& C \int^{\frac{M}{2}}_0(M-x)^{1-\frac{1}{\gamma}}W^2dx
      +C\int^M_{0}G\frac{M_0+x}{r_\infty^{2n-2}}W^2dx\nonumber\\
      &=&C \int^{\frac{M}{2}}_0(M-x)^{1-\frac{1}{\gamma}}W^2dx
      -C\int^M_0A(\rho_\infty^\gamma)_xW^2dx\nonumber\\
            &\leq&C \int^M_0\left((M-x)^\frac{1+\gamma}{\gamma}W_x^2
    +(M-x)^\frac{\gamma-1}{\gamma} W^2
    \right)dx,
    \end{eqnarray}
then, we have (\ref{ss-selfE2.8}) immediately.

Similarly, we obtain
 \begin{eqnarray*}
    &&S[V_\infty+W]-S[V_\infty] \nonumber\\
    &=&\frac{1}{2}\int^M_0\left\{ A[\gamma+O(|(M-x)^{\frac{1}
            {\gamma}}W_x|)]\rho_\infty^{1+\gamma}W_x^2
                -[2n-2+O(|W|)]G(M_0+x)r_\infty^{2-3n} W^2
    \right\}dx\nonumber\\
            &=&\frac{1}{2}\int^M_0\left\{ A[\gamma+O(|(M-x)^{\frac{1}
            {\gamma}}W_x|)]\rho_\infty^{1+\gamma}W_x^2
    \right.\nonumber\\
    &&\left.-2A[2n-2+O(|W|)]\rho_\infty^\gamma r_\infty^{-n} WW_x
    +nA[2n-2+O(|W|)]\rho_\infty^{\gamma-1}r_\infty^{-2n} W^2
    \right)dx,
    \end{eqnarray*}
for all  $W\in K_1$. Here, $O(d)$ means $O(d)\rightarrow 0$ as
$d\rightarrow0$. If $\gamma>\frac{2n-2}{n}$, choosing $\delta_9$
small enough, we can obtain (\ref{ss-self-E2.9}) immediately.
\end{proof}

Using the similar argument as that in Proposition
\ref{ss-self-stat-ex-prop3}, we could obtain the following
uniqueness result.
\begin{prop}\label{ss-self-stat-ex-prop3-1}
Let $\rho_\infty$ be a solution obtained in Proposition
\ref{ss-self-stat-ex-prop1}, and $\rho_2$ be another solution of
the problem (\ref{ss-E1.17})-(\ref{ss-E1.18}) satisfying
$\rho_2\in W^{1,\beta}([0,M])$ and
$C^{-1}(M-x)^\frac{1}{\gamma}\leq \rho_2(x)\leq
C(M-x)^\frac{1}{\gamma}$. If $\gamma>\frac{2n-2}{n}$ and
$\|(M-x)^{-\frac{1}{\gamma}}(\rho_\infty-\rho_2)(x)\|_{L^\infty}\leq
\delta_{10}$ with a small enough positive constant $\delta_{10}$,
then we have $\rho_\infty(x)=\rho_2(x)$, a.e. $x\in[0,M]$.
\end{prop}
\begin{proof}
From (\ref{ss-E1.17})-(\ref{ss-E1.18}), we have
    $$
    A\rho_\infty^\gamma=\int^M_x G\frac{M_0+y}{r_\infty^{2n-2}}dy,
      \ r_\infty^n(x)=a^n+n\int^x_0\rho_\infty^{-1}(y)dy,
    $$
    $$
    A\rho_2^\gamma=\int^M_x G\frac{M_0+y}{r_2^{2n-2}}dy,
      \ r_2^n(x)=a^n+n\int^x_0\rho_2^{-1}(y)dy,
    $$
and
    $$
    A(\rho_\infty^\gamma-\rho_2^\gamma)=\int^M_x G(M_0+y)\left(
    r_\infty^{2-2n}-r_2^{2-2n}\right)dy.
    $$
Multiplying the above equality by
$(\rho_\infty^{-1}-\rho_2^{-1})$, integrating  over $[0,M]$, and
using the fact that
$\int^M_0n(\rho_\infty^{-1}-\rho_2^{-1})(x)\int^M_x
g(y)dydx=\int^M_0 g(r_\infty^n-r_2^n)dx$, we obtain
    \begin{eqnarray*}
    0&=&\int^M_0\left\{A(\rho_\infty^\gamma-\rho_2^\gamma)(\rho_\infty^{-1}-\rho_2^{-1})
    -G\frac{(M_0+x)}{n}(r_\infty^{2-2n}-r_2^{2-2n})(r^n_\infty-r_2^n)
    \right\}dx\\
        &=&\int^M_0\left\{- A[\gamma+O(|(M-x)^{-\frac{1}
            {\gamma}}(\rho_\infty-\rho_2)|)]\rho_\infty^{1+\gamma}(\rho_\infty^{-1}-\rho_2^{-1})^2
            \right.\nonumber\\
        &&\left.-\frac{A}{n^2}[2n-2+O(|r_\infty^n-r_2^n|)](\rho_\infty^\gamma)_xr_\infty^{-n}(r_\infty^n-r_2^n)^2
    \right\}dx\nonumber\\
    &=&-\int^M_0\left\{ A[\gamma+O(|(M-x)^{-\frac{1}
            {\gamma}}(\rho_\infty-\rho_2)|)]\rho_\infty^{1+\gamma}(\rho_\infty^{-1}-\rho_2^{-1})^2
    \right.\nonumber\\
    &&-\frac{2A}{n}[2n-2+O(|r_\infty^n-r_2^n|)]\rho_\infty^\gamma r_\infty^{-n}
    (r_\infty^n-r_2^n)(\rho_\infty^{-1}-\rho_2^{-1})\\
    &&\left.+\frac{A}{n}[2n-2+O(|r_\infty^n-r_2^n|)]\rho_\infty^{\gamma-1}r_\infty^{-2n} (r_\infty^n-r^n_2)^2
    \right\}dx\\
            &\leq&-C^{-1}\int^M_0\left((M-x)^\frac{1+\gamma}{\gamma}(\rho_\infty^{-1}-\rho_2^{-1})^2
    +(M-x)^\frac{\gamma-1}{\gamma} (r_\infty^n-r_2^n)^2
    \right)dx,
    \end{eqnarray*}
when $\gamma>\frac{2n-2}{n}$ and $\delta_{10}$ is small enough.
Thus, we can obtain $\rho_\infty=\rho_2$ immediately.
\end{proof}

Now, we shall use the shooting method to prove the uniqueness of
the solution $\rho_\infty\in K$.
\begin{prop}\label{ss-self-stat-ex-prop4}
Under the assumption (\ref{ss-self-E2.1}), the Lagrangian
stationary problem (\ref{ss-E1.17})-(\ref{ss-E1.18}) has a unique
solution $\rho_\infty\in K$.
\end{prop}
\begin{proof}
 We
consider the Cauchy problem
    \begin{equation}
      (A\rho^\gamma)_x=-G(M_0+x)(nV)^{\frac{2-2n}{n}},
      \ (V)_x=\rho^{-1},
      \ x\in(0,M),\label{ss-self-E2.14}
    \end{equation}
        \begin{equation}
          \rho\big|_{x=0}=\sigma,
          \ V\big|_{x=0}=\frac{a^n}{n},\label{ss-self-E2.15}
        \end{equation}
for the unknown functions $\rho(\sigma,x)$ and $V(\sigma,x)$,
where $\sigma>0$ is the shooting parameter. Thus, for each
$\sigma>0$, using  the classical ODE theory, there exists a unique
solution to the problem
(\ref{ss-self-E2.14})-(\ref{ss-self-E2.15}) satisfying
$\rho(\sigma,x)>0$ for $x\in[0,M_\sigma)$, where either
$\rho\big|_{x=M_\sigma}=0$ and $M_\sigma\in(0,M)$ or $M_\sigma=M$.

Clearly, if $\rho_\infty\in K$ is a solution to the problem
(\ref{ss-E1.17})-(\ref{ss-E1.18}), then $\rho_\infty$ satisfies
(\ref{ss-self-E2.14})-(\ref{ss-self-E2.15}) for some $\sigma_0>0$,
and $M_{\sigma_0}=M$. We will show it is possible only for one
value of $\sigma$. Using  similar arguments as that in
\cite{Hartman} (\S V.3), we obtain that $(\partial_\sigma
\rho^\gamma,\partial_\sigma V)$ is  well defined and satisfies the
linear Cauchy problem
    \begin{equation}
      A(\partial_\sigma
\rho^\gamma)_x=(2n-2)G(M_0+x)(nV)^\frac{2-3n}{n}\partial_\sigma V,
\ (\partial_\sigma
V)_x=-\frac{1}{\gamma}\rho^{-\gamma-1}\partial_\sigma \rho^\gamma,
\ x\in[0,M_\sigma),\label{ss-self-E2.20}
    \end{equation}
        \begin{equation}
    \partial_\sigma
\rho^\gamma\big|_{x=0}=1,\ \partial_\sigma
V\big|_{x=0}=0.\label{ss-self-E2.21}
        \end{equation}
It is easy to see that
    $$
    \partial_\sigma
\rho^\gamma>0,\ (\partial_\sigma V)_x<0,
\partial_\sigma V<0
        $$
hold on $(0,M_4)$, where either $\partial_\sigma
\rho^\gamma\big|_{x=M_4}=0$ and $M_4\in(0,M_\sigma)$ or
$M_4=M_\sigma$. We claim that only $M_4=M_\sigma$ can occur.

Assume that $M_4\in(0,M_\sigma)$. Letting
$\phi=A\rho^\gamma(\partial_\sigma
V)_x+\frac{n}{2n-2}A\partial_\sigma \rho^\gamma(V)_x$, from
(\ref{ss-self-E2.14}) and (\ref{ss-self-E2.20}), we have
    $$
    \int^{M_4}_0\phi dx=\left.\left\{A\rho^\gamma\partial_\sigma
V+\frac{n}{2n-2}A\partial_\sigma \rho^\gamma
V\right\}\right|_0^{M_4}.
    $$
By the estimates $\rho(\sigma,M_4)>0$, $\partial_\sigma
\rho^\gamma\big|_{x=M_4}=0$, $\partial_\sigma V\big|_{x=M_4}<0$
and the initial condition (\ref{ss-self-E2.15}) and
(\ref{ss-self-E2.21}), we get
    $$
    \int^{M_4}_0\phi dx<0.
    $$
On the other hand, from (\ref{ss-self-E2.14}) and
(\ref{ss-self-E2.20}), we have
    $$
    \phi=A\rho^{-1}\partial_\sigma\rho^\gamma(\frac{n}{2n-2}-\frac{1}{\gamma})>0,
    \ x\in (0,M_4).
    $$
It is a contradiction.

Thus, we obtain
    $$
    \rho(\sigma,x)>0,\ \partial_\sigma(\sigma,x) \rho>0,
    \ x\in (0,M_\sigma),
    $$
and $M_\sigma$ is non-decreasing on $\sigma\in (0,\infty)$.
Therefore, for each fixed point $x\in [0,M_b)$, the function
$\rho(\sigma,x)$ is strictly increasing on $\sigma\geq b$.

If there exists $\sigma_1\not=\sigma_0$ such that
$M_{\sigma_1}=M_{\sigma_0}=M$ and $\rho(\sigma_1,x)\in K$, then
there exists
$\min\{\sigma_0,\sigma_1\}<\sigma_2<\max\{\sigma_0,\sigma_1\}$
such that
$0<\|(M-x)^{-\frac{1}{\gamma}}(\rho(\sigma_2,x)-\rho(\sigma_0,x))\|_{L^\infty}\leq
\delta_{10}$. From Proposition \ref{ss-self-stat-ex-prop3-1}, we
have $\rho(\sigma_2,x)=\rho(\sigma_0,x)=\rho_\infty(x)$, which is
a contradiction. Thus, we finish the proof of Proposition
\ref{ss-self-stat-ex-prop4}.
\end{proof}

Using the properties of the transformation (\ref{ss-self-E1.13})
and Propositions
\ref{ss-self-stat-ex-prop1}-\ref{ss-self-stat-ex-prop4}, we can
obtain the following proposition immediately.
\begin{prop}\label{ss-self-stat-ex-prop5}
  Under the assumption (\ref{ss-self-E2.1}),
the Eulerian stationary problem
(\ref{ss-E-sym1})-(\ref{ss-E-sym3}) has a unique  solution
$(\rho_\infty,l_\infty)$, satisfying $\rho_\infty(r)\sim
(l_\infty^n-r^n)^{\frac{1}{\gamma-1}}$, $(\rho_\infty)_r(r)<0$,
$a<r<l_\infty$ with $l_\infty<+\infty$.
\end{prop}
\begin{rem}
  The uniqueness of the solution in Proposition \ref{ss-self-stat-ex-prop5} means that:
  if $(\rho_{\infty1},l_{\infty1})$ and
  $(\rho_{\infty2},l_{\infty2})$ are two solutions to the Eulerian stationary problem
(\ref{ss-E-sym1})-(\ref{ss-E-sym3}) with the same total mass $M$,
and satisfy $\rho_{\infty i}(r)\sim (l_{\infty
i}^n-r^n)^{\frac{1}{\gamma-1}}$, $i=1,2$, then we
have$(\rho_{\infty1},l_{\infty1})=(\rho_{\infty2},l_{\infty2})$.
\end{rem}
\section{Global Existence}\label{ss-sec2}
Using  similar arguments as that in \cite{ChenP}, we  obtain the
following local existence and uniqueness result and omit the
proof.

\begin{thm}\label{ss-localthm}(Local Result) Under the assumptions in Theorem \ref{ss-thm},
there is a positive constant $T_{1}>0$ such that the free boundary
problem (\ref{ss-E})-(\ref{ss-Efd}) admits a unique weak solution
$(\rho, u, r)(x,t)$ on $[0, M]\times[0, T_{1}]$ in the sense that
$$\rho(x, t), u(x, t), r(x, t)\in L^\infty([0, M]\times[0, T_{1}])
  \cap C^{1}([0, T_{1}];L^{2}([0, M])),$$
$$\rho^{\theta+1}\partial_x(r^{n-1}u)\in L^{\infty}([0, M]\times[0, T_{1}])\cap
C^{\frac{1}{2}}([0, T_{1}];L^{2}([0, M])),$$
$$\partial_xr^{n-1}, \partial_x(r^{n-2}u) \in L^{\infty}([0, T_{1}], L^{1}([0, M])),$$
 and the following equations hold:
 $$\partial_{t}\rho=-\rho^{2}\partial_{x}(r^{n-1}u),\indent \rho(x,0)=\rho_0, $$
\begin{equation}
\partial_{t}r(x, t)=u(x, t),\ r^{n}(x, t)=a^{n}+n\int_{0}^{x}\rho^{-1}(y, t)dy,
\label{ss-2.2}
\end{equation}
    \begin{equation}
     (r^{\beta}(\rho^{\theta})_x)_t
    = -\frac{\theta r^{1+\beta-n}}{2c_{1}+c_{2}}u_t
    -\frac{\theta}{2c_{1}+c_{2}}
       \left(A
    r^{\beta}(\rho^{\gamma})_x
    +r^{1+\beta-n}f\right)
    \label{ss-2.3},
    \end{equation}
            \begin{eqnarray}
            &&(2c_1+c_2)\rho^{1+\theta}(r^{n-1}u)_x\nonumber\\
                    &=&A\rho^\gamma+2c_1(n-1)\rho^\theta\frac{u}{r}
            +\int^M_x\left\{-\frac{u_t}{r^{n-1}}
                +2c_1(n-1)\rho^\theta\left(\frac{u}{r}\right)_x-\frac{f}{r^{n-1}}\right\}dy
                 \label{ss-2.4},
            \end{eqnarray}
 for almost all $x\in[0, M
 ]$,  any $t\in[0,T_1]$, where
 $\beta=\frac{2(n-1)c_1\theta}{2c_1+c_2}$,
 \begin{eqnarray}
    &&\int^\infty_0\int^M_0[u\psi_t+(P-\rho(\lambda+2\mu)(r^{n-1}
    u)_x)(r^{n-1}\psi)_x\nonumber\\
    &&+2(n-1)\mu(r^{n-2}u\psi)_x-f(x,r,t)\psi]dxdt
        +\int^M_0u_0(x)\psi(x,0)dx=0,
        \end{eqnarray}
for any test function $\psi(x,t)\in
C^\infty_0((0,M]\times[0,T_1))$. Furthermore, we have
        \begin{equation}
\frac{N_1}{3}(1-x)^{\frac{1}{\gamma}}\leq\rho(x, t)
 \leq3N_2(1-x)^{\frac{1}{\gamma}} ,\  (x,t)\in[0,M]\times[0,T_1],
 \label{ss-E3.8}
       \end{equation}
    \begin{equation}
    (M-x)^{-\frac{1}{\gamma}}\rho(x,
    t)\in C([0,T_1];L^\infty([0,M])),\label{ss-E3.5}
    \end{equation}
    $$
    (M-x)^{\frac{\gamma-\theta}{2\gamma}}(\rho^\theta)_x,
    \ \rho_t,\ u_t\in L^\infty([0,T_1];L^2([0,M])),
    $$
    $$
      \rho^\frac{\theta+1}{2}u_{xt}\in  L^2([0,M]\times[0,T_1]),
     \  \rho\partial_x u \in
      L^\infty([0,M]\times[0,T_1]).
    $$
\end{thm}

\begin{rem}
  From (\ref{ss-E})$_1$, (\ref{ss-E3.8}) and $\rho\partial_x u \in
      L^\infty_{t,x}$, we have $(M-x)^{-\frac{1}{\gamma}}\partial_t\rho\in
      L^\infty_{tx}$. Thus, (\ref{ss-E3.5}) holds.
\end{rem}

Assume  the maximum existence time of the weak solution in Theorem
\ref{ss-localthm} is $T^*$. In this section, under the small
assumptions on the initial data, we will obtain the following
\textit{a priori} estimates and prove that $T^*=\infty$. In the
following, we may assume that $(\rho,u,r)(x,t)$ is suitably smooth
since the following estimates are valid for the solutions with the
regularities indicated in Theorem \ref{ss-localthm} by using the
Friedrichs mollifier.

From (\ref{ss-f2}), (\ref{ss-E1.17}) and Proposition
\ref{ss-self-stat-ex-prop1}, we could obtain the following lemma
easily.
\begin{lem}\label{ss-L2.1}
Under the assumptions of Theorem \ref{ss-thm}, we have
    \begin{equation}
      A\rho_\infty^\gamma(x)=\int^M_xG\frac{M_0+y}{r_\infty^{2n-2}}dy,
      \label{ss-rhoinf1}
    \end{equation}
    \begin{equation}
  C^{-1}(M-x)^\frac{1}{\gamma}\leq \rho_\infty\leq C(M-x)^\frac{1}{\gamma},
  \ r_\infty(x)\in[a,C], \label{ss-rhoinf}
  \end{equation}
        \begin{equation}
      \frac{d}{dx}\left(A\rho_\infty^\gamma(x)
      \right)=-G\frac{M_0+x}{r_\infty^{2n-2}},
      \ C^{-1}\leq (M-x)^{1-\frac{1}{\gamma}}\frac{d}{dx}\rho_\infty(x)\leq C,\label{ss-rhoinf3}
    \end{equation}
for all $x\in[0,M]$.
\end{lem}

\begin{lem}
Under the assumptions of Theorem \ref{ss-thm}, we have
    \begin{eqnarray}
      &&\frac{d}{dt}\int^M_0\left(\frac{1}{2}u^2+\frac{A\rho^{\gamma-1}}{\gamma-1}+
      \int^r_aG\frac{M_0+x}{s^{n-1}}ds\right)dx\nonumber\\
              &&+\int^M_0\left\{
              \left(\frac{2}{n}c_1+c_2\right)\rho^{1+\theta}[\partial_x(r^{n-1}u)]^2+\frac{2(n-1)}{n}
              c_1\rho^{1+\theta}(r^{n-1}u_x-\frac{u}{r\rho})^2\right\}dx\nonumber\\
      &=&-\int^M_0\Delta f
      udx, \ t\in [0,T^*).
       \label{ss-E2.1-0}
    \end{eqnarray}
\end{lem}
\begin{proof}
  Multiplying (\ref{ss-E})$_2$ by $u$, integrating the
resulting equation over $[0,M]$, using integration by parts and
 the boundary conditions
(\ref{ss-Efd}), we obtain
    \begin{eqnarray}
      &&\frac{d}{dt}\int^M_0\frac{1}{2}u^2dx
      +\int^M_0\left\{
      (2c_1+c_2)\rho^{1+\theta}[\partial_x(r^{n-1}u)]^2-
      2c_1(n-1)\rho^\theta\partial_x(r^{n-2}u^2)\right\}dx\nonumber\\
            &=&\int^M_0
            A\rho^\gamma\partial_x(r^{n-1}u)dx
            -\int^M_0     fudx.
            \label{ss-E2.2}
    \end{eqnarray}
From (\ref{ss-E}), we have
    \begin{equation}
    \int^M_0            A\rho^\gamma\partial_x(r^{n-1}u)dx
            =-\frac{d}{dt}\int^M_0\frac{A}{\gamma-1}\rho^{\gamma-1}dx,
    \end{equation}
            \begin{equation}
      -\int^M_0fudx=-\frac{d}{dt}\int^M_0\int^r_aG\frac{M_0+x}{s^{n-1}}dsdx-\int^M_0\Delta
      fudx,
    \end{equation}
and
    \begin{eqnarray}
      &&
      (2c_1+c_2)\rho^{1+\theta}[\partial_x(r^{n-1}u)]^2-
      2c_1(n-1)\rho^\theta\partial_x(r^{n-2}u^2)\nonumber\\
            &=&
     \left(\frac{2}{n}c_1+c_2\right)\rho^{1+\theta}[\partial_x(r^{n-1}u)]^2+\frac{2(n-1)}{n}
      c_1\rho^{1+\theta}(r^{n-1}u_x-\frac{u}{r\rho})^2.\label{ss-E2.5}
    \end{eqnarray}
From (\ref{ss-E2.2})-(\ref{ss-E2.5}), we obtain
 (\ref{ss-E2.1-0}) immediately.
\end{proof}

Now, using the classical continuity method, we will obtain the
estimate of
$\|(M-x)^{-\frac{1}{\gamma}}(\rho-\rho_\infty)\|_{L^\infty}$.

\noindent\textbf{Claim 1}: Under the assumptions of Theorem
\ref{ss-thm}, there is a small positive constant
$\epsilon_1>\epsilon_0$, such that, for any $T\in (0,T^*)$, if
    \begin{equation}
      I(t)=\|g(\cdot,t)-g_\infty\|_{L^\infty}     \leq
      2\epsilon_1,\ \forall\ t\in[0,T],\label{ss-E2.7}
    \end{equation}
where $g(x,t)=(M-x)^{-\frac{1}{\gamma}}\rho(x,t)$ and
$g_\infty(x)=(M-x)^{-\frac{1}{\gamma}}\rho_\infty(x)$, then
    \begin{equation}
    I(t)\leq \epsilon_1, \ \forall\ t\in[0,T].
    \end{equation}

    Using the results in Lemmas \ref{ss-L2.2}-\ref{ss-Lem3.4}, we
    can give the definition of $\epsilon_1$ in
    (\ref{ss-self-E3.76}) and finish the proof of \textbf{Claim
    1}.
\begin{lem}\label{ss-L2.2}
  Under the assumptions of Theorem \ref{ss-thm} and
  (\ref{ss-E2.7}), if $\epsilon_1$ is small enough, we obtain
        \begin{equation}
          C_1^{-1}(M-x)^\frac{1}{\gamma}\leq
          \rho(x,t)\leq           C_1(M-x)^\frac{1}{\gamma},
          \label{ss-E2.9}
        \end{equation}
  \begin{equation}
    r(x,t)\in [a,C_1],
              \label{ss-E2.10}
  \end{equation}
for all $t\in[0,T]$ and $x\in[0,M]$.
\end{lem}
\begin{proof}
    From (\ref{ss-E})$_3$, (\ref{ss-E2.7}) and Lemma \ref{ss-L2.1},
   we can easily obtain the estimate (\ref{ss-E2.9})  and
(\ref{ss-E2.10}) when $4\epsilon_1<\min_{x\in[0,M]}g_\infty$.
\end{proof}
\begin{lem}
  Under the assumptions of Lemma \ref{ss-L2.2}, if $\epsilon_1$ is small enough, we obtain
        \begin{equation}
      \int^M_0\left\{u^2+(M-x)^{\frac{\gamma-1}{\gamma}}(g-g_\infty)^2
      +(r-r_\infty)^2\right\}dx\leq C_2\epsilon_0^2,
                    \label{ss-E2.11}
        \end{equation}
        \begin{equation}
    \int^t_0\|u(\cdot,s)\|_{L^\infty}^2ds+    \int^t_0\int^M_0\left(\rho^{\theta+1}u_x^2+\rho^{\theta-1}u^2
        \right)(x,s)dxds\leq C_2\epsilon_0^2,
        \label{ss-E2.13}
        \end{equation}
for all $t\in[0,T]$.
\end{lem}

\begin{proof}
 From
(\ref{ss-self-E2.9-1}), (\ref{ss-rhoinf1}) and (\ref{ss-E2.1-0}),
we have
    \begin{eqnarray}
      &&\frac{d}{dt}\left(\int^M_0\frac{1}{2}u^2dx+S[V]-S[V_\infty]\right)\nonumber\\
              &&+\int^M_0\left\{
              \left(\frac{2}{n}c_1+c_2\right)\rho^{1+\theta}[\partial_x(r^{n-1}u)]^2+\frac{2(n-1)}{n}
              c_1\rho^{1+\theta}(r^{n-1}u_x-\frac{u}{r\rho})^2\right\}dx\nonumber\\
      &=&
    -\int^M_0\Delta f udx
     \label{ss-E2.16-1}
    \end{eqnarray}
where $V_\infty=\frac{r_\infty^n}{n}$ and $V=\frac{r^n}{n}$. From
(\ref{ss-self-E2.9}), (\ref{ss-E2.9})-(\ref{ss-E2.10}) and
Proposition \ref{ss-self-stat-ex-prop3}, we have
        \begin{eqnarray}
      &&C^{-1}\int^M_0(M-x)^{\frac{\gamma-1}{\gamma}}(g-g_\infty)^2
      +(r-r_\infty)^2dx\nonumber\\
      &\leq&  S[V]-S[V_\infty]\leq
      C\int^M_0(M-x)^{\frac{\gamma-1}{\gamma}}(g-g_\infty)^2
      +(r-r_\infty)^2dx,\label{ss-E2.18-1}
        \end{eqnarray}
when
$\|(M-x)^{\frac{1}{\gamma}}(\rho^{-1}-\rho_\infty^{-1})\|_{L^\infty_x}
+\|\frac{1}{n}(r^n-r_\infty^n)\|_{L^\infty_x}\leq
C_3\epsilon_1\leq\delta_9$. From (\ref{ss-epsilon2}),
(\ref{ss-E2.9})-(\ref{ss-E2.10}) and
(\ref{ss-E2.16-1})-(\ref{ss-E2.18-1}), we obtain
    \begin{eqnarray}
      &&\int^M_0\left\{u^2+(M-x)^{\frac{\gamma-1}{\gamma}}(g-g_\infty)^2
      +(r-r_\infty)^2\right\}dx+\int^t_0\int^M_0\left\{\rho^{\theta+1}u_x^2
              +\rho^{\theta-1}u^2\right\}dxds\nonumber\\
      &\leq&C\epsilon_0^2+C\int^t_0
      f_1(s)\|u(\cdot,s)\|_{L^\infty}ds.\label{ss-E3.22}
    \end{eqnarray}
Since $\theta\in(0,\gamma-1)$, we obtain
    \begin{eqnarray}
      |u(x,t)|&=&\left|\int^x_0u_xdy
      \right|\leq C\left(\int^x_0\rho^{\theta+1}u_x^2dy
      \right)^\frac{1}{2}\left(\int^x_0\rho^{-\theta-1}dy
      \right)^\frac{1}{2}\nonumber\\
        &\leq& C\left(\int^x_0\rho^{\theta+1}u_x^2dy
      \right)^\frac{1}{2}\left(\int^x_0(M-y)^{-\frac{\theta+1}{\gamma}}dy
      \right)^\frac{1}{2}
      \leq C\left(\int^x_0\rho^{\theta+1}u_x^2dy
      \right)^\frac{1}{2}\label{ss-E3.32}
    \end{eqnarray}
and
    \begin{equation}
    C\int^t_0
      f_1(s)\|u(\cdot,s)\|_{L^\infty}ds\leq
      \frac{1}{2}\int^t_0\int^M_0\rho^{\theta+1}u_x^2dyds
      +C^2\int^t_0f_1^2dt.\label{ss-E3.23}
    \end{equation}
From (\ref{ss-epsilon1}) and (\ref{ss-E3.22})-(\ref{ss-E3.23}), we
can obtain (\ref{ss-E2.11})-(\ref{ss-E2.13}) immediately.
\end{proof}
\begin{lem}\label{ss-L3.5-0}
  Under the assumptions of Lemma \ref{ss-L2.2}, if $\epsilon_1$ is small enough, we obtain
        \begin{equation}
   (1+t)^\alpha\int^M_0\rho_\infty^{\theta-1}
    (g-g_\infty)^2dx+ \int^{t}_0\int^M_0
     (1+s)^\alpha \left[ \rho_\infty^{\gamma-1}
    (g-g_\infty)^2+(r-r_\infty)^2\right]dxds\leq
    C_4\epsilon_0,\label{ss-E2.21-2}
        \end{equation}
for all $t\in[0,T]$, where $\alpha=-\frac{5}{8}$.
\end{lem}
\begin{proof}
  Multiplying (\ref{ss-E})$_2$ by $(1+t)^\alpha r^{1-n}(\frac{r^n}{n}-\frac{r^n_\infty}{n})$, integrating
   over $[0,M]$, using integration by parts and
 the boundary conditions
(\ref{ss-Efd}), we obtain
    \begin{eqnarray}
      &&(1+t)^\alpha\int^M_0\left[A(\rho_\infty^\gamma-\rho^\gamma)(\rho^{-1}-\rho_\infty^{-1})
      +G(M_0+x)(r^{2-2n}-r_\infty^{2-2n})(\frac{r^n}{n}-\frac{r^n_\infty}{n})\right]dx\nonumber\\
            &=&-(1+t)^\alpha\int^M_0\frac{u_t}{r^{n-1}}(\frac{r^n}{n}-\frac{r^n_\infty}{n})dx
            -(1+t)^\alpha\int^M_0\Delta fr^{1-n}\left(\frac{r^n}{n}-\frac{r^n_\infty}{n}
      \right)dx
            \nonumber\\
            &&+(1+t)^\alpha\int^M_0
            (2c_1+c_2)\rho^{1+\theta}\partial_x(r^{n-1}u)(\rho_\infty^{-1}-\rho^{-1})dx\nonumber\\
     &&+(1+t)^\alpha\int^M_0
            2c_1(n-1)\rho^\theta \left(\frac{u}{r}(\frac{r^n}{n}-\frac{r^n_\infty}{n})
                   \right)_xdx:=\sum^4_{i=1}B_i.\label{ss-E2.22-1}
    \end{eqnarray}
We can rewrite the left hand side of (\ref{ss-E2.22-1}) as follows
    \begin{eqnarray*}
      \textrm{L.H.S of } (\ref{ss-E2.22-1})
      &=&(1+t)^\alpha\int^M_0
      \left[
      A(\gamma+O(\epsilon_1))\rho_\infty^{\gamma+1}(\rho^{-1}-\rho_\infty^{-1})^2\right.\\
      &&\left.-      (2n-2+O(\epsilon_1))G(M_0+x)r_\infty^{2-3n}
      (\frac{r^n}{n}-\frac{r_\infty^n}{n})^2
      \right]dx.
    \end{eqnarray*}
Similar to (\ref{ss-selfE2.8}), we have
    \begin{equation}
    \textrm{L.H.S of } (\ref{ss-E2.22-1})\geq C_5(1+t)^\alpha\int^M_0
      \left[\rho_\infty^{\gamma-1}(g-g_\infty)^2+
      \left(r-r_\infty\right)^2
      \right]dx,\label{ss-E2.22-1-2}
    \end{equation}
when $\epsilon_1\leq \delta_{10}$ is small enough.

Using  (\ref{ss-E2.9})-(\ref{ss-E2.11}),
 integration by parts and H\"{o}lder's inequality, we can
estimate $B_i$ as follows.
\begin{eqnarray}
  B_1&=&-\frac{d}{dt}\int^M_0(1+t)^\alpha\frac{u}{r^{n-1}}\left(\frac{r^n}{n}-\frac{r_\infty^n}{n}
      \right)dx+\alpha(1+t)^{\alpha-1}\int^M_0\frac{u}{r^{n-1}}\left(\frac{r^n}{n}-\frac{r_\infty^n}{n}
      \right)dx\nonumber\\
      &&+(1+t)^\alpha\int^M_0u^2\left(\frac{1}{n}+\frac{(n-1)r_\infty^n}{nr^n}
      \right)dx\nonumber\\
            &\leq& -\frac{d}{dt}\int^M_0(1+t)^\alpha\frac{u}{r^{n-1}}\left(\frac{r^n}{n}-\frac{r_\infty^n}{n}
      \right)dx+
      C\int^M_0u^2dx+C\epsilon_0^2(1+t)^{\alpha-1},\label{ss-E3.29}
\end{eqnarray}
    \begin{equation}
      B_2\leq C\epsilon_0(1+t)^\alpha f_1,
    \end{equation}
    \begin{eqnarray}
      B_3&=&-\frac{2c_1+c_2}{\theta}\int^M_0(\rho^\theta)_t(1+t)^\alpha
      (\frac{1}{\rho_\infty}-\frac{1}{\rho})dy\nonumber\\
      &=&-\frac{2c_1+c_2}{\theta}\int^M_0\partial_th(\rho,\rho_\infty)(1+t)^\alpha
      dx\nonumber\\
            &=&-\frac{2c_1+c_2}{\theta}\frac{d}{dt}\int^M_0h(\rho,\rho_\infty)(1+t)^\alpha
      dx+\frac{\alpha(2c_1+c_2)}{\theta}\int^M_0h(\rho,\rho_\infty)(1+t)^{\alpha-1}
      dx,
    \end{eqnarray}
where $h(\rho,\rho_\infty)=\int^\rho_{\rho_\infty}\theta
s^{\theta-1}(\frac{1}{\rho_\infty}-\frac{1}{s})ds\sim
\rho_\infty^{\theta-1}(g-g_\infty)^2$, and
    \begin{eqnarray}
      B_4&\leq&C(1+t)^{\alpha}\int^M_0
      \left[\rho^{\theta}|u_x|+\rho^{\theta-1}|u|\right]dx\nonumber\\
                &\leq&       C(1+t)^{\alpha}\left[\int^M_0
                \left(\rho^{\theta+1}u_x^2+\rho^{\theta-1}u^2
                \right)dx      \right]^\frac{1}{2}
                \left[\int^M_0(M-x)^{\frac{\theta-1}{\gamma}}
                \right]^{\frac{1}{2}}\nonumber\\
    &\leq&       C(1+t)^{\alpha}\left[\int^M_0
                \left(\rho^{\theta+1}u_x^2+\rho^{\theta-1}u^2
                \right)dx      \right]^\frac{1}{2},\label{ss-E2.30-1}
    \end{eqnarray}
since $\gamma+\theta-1>0$. From
(\ref{ss-E2.22-1})-(\ref{ss-E2.30-1}), we get
    \begin{eqnarray*}
    &&\frac{d}{dt}\int^M_0(1+t)^\alpha\left\{\frac{u}{r^{n-1}}\left(\frac{r^n}{n}-\frac{r_\infty^n}{n}
      \right)+\frac{2c_1+c_2}{\theta}h(\rho,\rho_\infty)
      \right\}dx\nonumber\\
      &&+C^{-1}\int^M_0(1+t)^{\alpha-1}\rho_\infty^{\theta-1}(g-g_\infty)^2+(1+t)^\alpha\left\{\rho_\infty^{\gamma-1}\left[
    (g-g_\infty)^2+(r-r_\infty)^2\right]\right\}dx\nonumber\\
      &  \leq& C(1+t)^{\alpha}\left[\int^M_0\rho_\infty^{\theta+1}u_x^2dx+\|u(\cdot,t)\|_{L^\infty}^2
      \right]^\frac{1}{2}+
      C\epsilon_0(1+t)^\alpha f_1+C\int^M_0u^2dx+C\epsilon_0^2(1+t)^{\alpha-1}          .
    \end{eqnarray*}
And using (\ref{ss-epsilon1}), (\ref{ss-E2.9})-(\ref{ss-E2.13})
and H\"{o}lder's inequality,  we obtain (\ref{ss-E2.21-2})
immediately.
\end{proof}

Let $\epsilon_2\in(0,\min\{\frac{1}{4}$,
$\frac{\gamma-\theta-1}{\gamma-\theta},\frac{\gamma-1}{2(3\gamma-1)}\})$
be a constant. Define $\{\beta_j\}$ and $\{\alpha_j\}$ by
$\beta_{j+1}=\frac{\beta_j}{2}+\frac{1}{2}-\frac{\epsilon_2}{4}$,
$\alpha_j=\frac{\beta_j}{2}-\frac{1}{2}-\frac{\epsilon_2}{4}$ and
$\alpha_0=\alpha=-\frac{5}{8}$, $j=0,1,\ldots$. Let $N_4$ be a
integer satisfying
$\beta_{N_4}\in[1-\epsilon_2,1-\frac{3\epsilon_2}{4})$ and
$\alpha_{N_4}\in(-\epsilon_2,-\frac{\epsilon_2}{4})$. It is easy
to see that $\beta_0=-\frac{1}{4}+\frac{\epsilon_2}{2}<0$,
$\alpha_j\in[-\frac{5}{8},-\frac{\epsilon_2}{4})$ and
$\beta_j\in(-\frac{1}{4},1-\frac{3\epsilon_2}{4})$,
$j=0,1,\ldots,N_4$. Then, the following lemma can be proved by
induction.

\begin{lem}\label{ss-L3.7}
 Under the assumptions of Lemma \ref{ss-L2.2}, if $\epsilon_1$ is small enough,  we obtain
  \begin{equation}
      \int^M_0\left\{u^2+(M-x)^{\frac{\gamma-1}{\gamma}}(g-g_\infty)^2
      +(r-r_\infty)^2\right\}dx\leq
      C_7\epsilon_0^{2^{1-N_4}}(1+t)^{\epsilon_2-1},\label{ss-E4.1-0}
        \end{equation}
        \begin{equation}
    \int^t_0(1+s)^{1-\epsilon_2}\|u(\cdot,s)\|_{L^\infty}^2ds+    \int^t_0\int^M_0
    (1+s)^{1-\epsilon_2}\left(\rho^{\theta+1}u_x^2+\rho^{\theta-1}u^2
        \right)(x,s)dxds\leq C_7\epsilon_0^{2^{1-N_4}},\label{ss-E4.2-0}
        \end{equation}
            \begin{equation}
   \int^M_0\frac{(M-x)^\frac{\theta-1}{\gamma}}{(1+t)^{\epsilon_2}}(g-g_\infty)^2dx+ \int^{t}_0\int^M_0
    \frac{\rho_\infty^{\gamma-1}(g-g_\infty)^2+(r-r_\infty)^2}{(1+s)^{\epsilon_2} }dxds\leq
    C_7\epsilon_0^{2^{-N_4}},\label{ss-E4.2-3}
        \end{equation}
 for all $t\in[0,T]$, where $C_7$ is a constant depending on $\epsilon_2$.
\end{lem}
\begin{proof}The following
estimates can be proved by induction:
 \begin{equation}
      \int^M_0\left\{u^2+(M-x)^{\frac{\gamma-1}{\gamma}}(g-g_\infty)^2
      +(r-r_\infty)^2\right\}dx\leq
      C_{j,\epsilon_2}\epsilon_0^{2^{1-j}}(1+t)^{-\beta_j},\label{ss-E4.1}
        \end{equation}
        \begin{equation}
    \int^t_0(1+s)^{\beta_j}\|u(\cdot,s)\|_{L^\infty}^2ds+    \int^t_0\int^M_0
    (1+s)^{\beta_j}\left(\rho^{\theta+1}u_x^2+\rho^{\theta-1}u^2
        \right)(x,s)dxds\leq C_{j,\epsilon_2}\epsilon_0^{2^{1-j}},\label{ss-E4.2}
        \end{equation}
            \begin{equation}
   (1+t)^{\alpha_j}\int^M_0\rho_\infty^{\theta-1}(g-g_\infty)^2dx+ \int^{t}_0\int^M_0
    (1+s)^{\alpha_j}\left[ \rho_\infty^{\gamma-1}
    (g-g_\infty)^2+(r-r_\infty)^2\right]dxds\leq
    C_{j,\epsilon_2}\epsilon_0^{2^{-j}},\label{ss-E4.3}
        \end{equation}
for all $t\geq 0$, where $C_{j,\epsilon_2}$ is a constant
depending on $j$ and $\epsilon_2$,
 $j=0,1,\ldots,N_4$.

  From (\ref{ss-E2.11})-(\ref{ss-E2.13}) and
 (\ref{ss-E2.21-2}), we obtain
 (\ref{ss-E4.1})-(\ref{ss-E4.3}) hold with $j=0$. Now, suppose
 that (\ref{ss-E4.1})-(\ref{ss-E4.3}) hold with $j=k\geq0$. To
 show (\ref{ss-E4.1})-(\ref{ss-E4.2}) hold with $j=k+1$, from (\ref{ss-E2.16-1}), we have
    \begin{eqnarray*}
      &&\frac{d}{dt}\left\{(1+t)^{\beta_{k+1}}\left(\int^M_0\frac{1}{2}u^2dx+S[V]-S[V_\infty]\right)
      \right\}\nonumber\\
              &&+(1+t)^{\beta_{k+1}}\int^M_0\left\{
              (\frac{2}{n}c_1+c_2)\rho^{1+\theta}[\partial_x(r^{n-1}u)]^2+\frac{2(n-1)}{n}
              c_1\rho^{1+\theta}(r^{n-1}u_x-\frac{u}{r\rho})^2\right\}dx\nonumber\\
      &=&\beta_{k+1}(1+t)^{\alpha_k}\left(\int^M_0\frac{1}{2}u^2(x,t)dx+S[V]-S[V_\infty]\right)
    -(1+t)^{\beta_{k+1}}\int^M_0\Delta f udx,
    \end{eqnarray*}
where $V_\infty=\frac{r_\infty^n}{n}$ and $V=\frac{r^n}{n}$.
Integrating the above equality in $[0,t]$, using
(\ref{ss-epsilon2}), (\ref{ss-E2.9})-(\ref{ss-E2.13}),
 (\ref{ss-E2.18-1}),
(\ref{ss-E4.3}) with $j=k$ and the fact that $\alpha_{k}<0$, we
obtain
    \begin{eqnarray}
      &&(1+t)^{\beta_{k+1}}\int^M_0\left\{u^2(x,t)+(M-x)^{\frac{\gamma-1}{\gamma}}(g-g_\infty)^2
      +(r-r_\infty)^2\right\}dx\nonumber\\
              &&+\int^t_0\int^M_0(1+s)^{\beta_{k+1}}\left\{\rho^{\theta+1}u_x^2
              +\rho^{\theta-1}u^2\right\}dxds\nonumber\\
      &\leq&C\epsilon_0^{2^{-k}}+C\int^t_0
      (1+s)^{\beta_{k+1}}f_1(s)\|u(\cdot,s)\|_{L^\infty}ds.\label{ss-E3.51}
    \end{eqnarray}
From (\ref{ss-epsilon1}) and (\ref{ss-E3.23}), we can obtain
(\ref{ss-E4.1})-(\ref{ss-E4.2}) with $j=k+1$ immediately.

To show (\ref{ss-E4.3}) with $j=k+1$, from
(\ref{ss-E2.22-1})-(\ref{ss-E2.22-1-2}), we have
 \begin{eqnarray}
      &&(1+t)^{\alpha_{k+1}}\int^M_0
      \left[A(\rho_\infty^\gamma-\rho^\gamma)(\rho^{-1}-\rho_\infty^{-1})
      +G(M_0+x)(r^{2-2n}-r_\infty^{2-2n})(\frac{r^n}{n}-\frac{r^n_\infty}{n})\right]dx\nonumber\\
            &=&-(1+t)^{\alpha_{k+1}}\int^M_0\frac{u_t}{r^{n-1}}(\frac{r^n}{n}-\frac{r^n_\infty}{n})dx
            -(1+t)^{\alpha_{k+1}}\int^M_0\Delta fr^{1-n}\left(\frac{r^n}{n}-\frac{r^n_\infty}{n}
      \right)dx
            \nonumber\\
            &&+(1+t)^{\alpha_{k+1}}\int^M_0
            (2c_1+c_2)\rho^{1+\theta}\partial_x(r^{n-1}u)(\rho_\infty^{-1}-\rho^{-1})dx\nonumber\\
     &&+(1+t)^{\alpha_{k+1}}\int^M_0
            2c_1(n-1)\rho^\theta \left(\frac{u}{r}(\frac{r^n}{n}-\frac{r^n_\infty}{n})
                   \right)_xdx:=\sum^4_{i=1}E_i,\label{ss-E2.22-5}
    \end{eqnarray}
and
    \begin{equation}
    \textrm{L.H.S of } (\ref{ss-E2.22-5})\geq C_{8}(1+t)^{\alpha_{k+1}}\int^M_0
     \left[ \rho_\infty^{\gamma-1}(g-g_\infty)^2+
      \left(r-r_\infty\right)^2
      \right]dx.
    \end{equation}
Similar to (\ref{ss-E3.29})-(\ref{ss-E2.30-1}), applying the
estimates (\ref{ss-E2.9})-(\ref{ss-E2.11}), integration by parts,
H\"{o}lder's inequality and the fact that $\alpha_{k+1}<0$, we can
estimate $E_i$ as follows.
\begin{equation}
  E_1\leq -\frac{d}{dt}\int^M_0(1+t)^{\alpha_{k+1}}\frac{u}{r^{n-1}}\left(\frac{r^n}{n}-\frac{r_\infty^n}{n}
      \right)dx+ C\|u\|_{L^\infty_x}^2+C\epsilon_0^2(1+t)^{{\alpha_{k+1}}-1},
\end{equation}
    \begin{equation}
      E_2\leq C\epsilon_0 f_1(1+t)^{\alpha_{k+1}},
    \end{equation}
    \begin{equation}
      E_3
            =-\frac{2c_1+c_2}{\theta}\frac{d}{dt}\int^M_0h(\rho,\rho_\infty)(1+t)^{\alpha_{k+1}}
      dx+\frac{{\alpha_{k+1}}(2c_1+c_2)}{\theta}\int^M_0h(\rho,\rho_\infty)(1+t)^{{\alpha_{k+1}}-1}
      dx,
    \end{equation} and
    \begin{equation}
      E_4\leq C(1+t)^{-\frac{1}{2}-\frac{\epsilon_2}{4}}\left[(1+t)^{\beta_{k+1}}
      \int^M_0\rho^{\theta+1}u_x^2dx+\|u(\cdot,t)\|_{L^\infty}^2\right]^\frac{1}{2}.\label{ss-E2.30-7}
    \end{equation}
Using (\ref{ss-epsilon1}), (\ref{ss-E2.22-5})-(\ref{ss-E2.30-7}),
 (\ref{ss-E4.1})-(\ref{ss-E4.2}) with $j=k+1$ and H\"{o}lder's
inequality, we get
    \begin{eqnarray}
    &&\int^M_0(1+t)^{\alpha_{k+1}}\rho_\infty^{\theta-1}(g-g_\infty)^2
      dx+\int^t_0\int^M_0(1+s)^{{\alpha_{k+1}}-1}\rho_\infty^{\theta-1}(g-g_\infty)^2dxds\nonumber\\
    &&+\int^t_0\int^M_0
      (1+s)^{\alpha_{k+1}}\left\{\rho_\infty^{\gamma-1}\left[
    (g-g_\infty)^2+(r-r_\infty)^2\right]\right\}dxds\nonumber\\
      &  \leq& C(1+t)^{\alpha_{k+1}}\int^M_0|u||r-r_\infty|dx+C\int^M_0|u_0||r_0-r_\infty|dx\nonumber\\
        &&+C\int^t_0\left[f_1(1+s)^{\alpha_{k+1}}+\epsilon_0^2(1+s)^{{\alpha_{k+1}}-1}+
        \|u\|_{L^\infty_x}^2\right]ds\nonumber\\
        &&+ C\int^t_0(1+s)^{-\frac{1}{2}-\frac{\epsilon_2}{4}}\left[
     (1+t)^{\beta_{k+1}}\left( \int^M_0\rho^{\theta+1}u_x^2dx+\|u(\cdot,t)\|_{L^\infty}^2
     \right)\right]^\frac{1}{2}ds\nonumber\\
        &\leq&C \epsilon_0^{2^{-(k+1)}},
          \label{ss-E2.33-7}
    \end{eqnarray}
and finish the proof of (\ref{ss-E4.3}) with $j=k+1$. Thus, we
show that (\ref{ss-E4.1})-(\ref{ss-E4.3}) hold for $j=0,1,\ldots,
N_4$, and obtain (\ref{ss-E4.1-0})-(\ref{ss-E4.2-3}) immediately.
\end{proof}

From Lemma \ref{ss-L3.7}, we can obtain the following estimate of
the weighted $L^2-$norm of $g-g_\infty$.

\begin{lem}\label{ss-L3.8}
 Under the assumptions of Lemma \ref{ss-L2.2}, we obtain
    \begin{equation}
    \int^M_0(M-x)^\frac{\theta-1+(\gamma-\theta)\epsilon_2}{\gamma}(g-g_\infty)^2dx
    \leq C_9\epsilon_0^{2^{-N_4}}, \ t\in[0,T].\label{ss-E3.66}
    \end{equation}
\end{lem}
\begin{proof}
Using (\ref{ss-E4.1-0}), (\ref{ss-E4.2-3}) and H\"{o}lder's
inequality, we have
    \begin{eqnarray*}
        && \int^M_0(M-x)^\frac{\theta-1+(\gamma-\theta)\epsilon_2}{\gamma}(g-g_\infty)^2dx\nonumber\\
            &\leq& C\left[\int^M_0
            (1+t)^{1-\epsilon_2}(M-x)^{\frac{\gamma-1}{\gamma}}(g-g_\infty)^2dx
            \right]^{\epsilon_2}
            \left[\int^M_0\frac{(M-x)^\frac{\theta-1}{\gamma}}{(1+t)^{\epsilon_2}}(g-g_\infty)^2dx
            \right]^{1-\epsilon_2}\nonumber\\
                &\leq& C\epsilon_0^{2^{-N_4}}.
    \end{eqnarray*}
\end{proof}

Then, using the similar argument as that in \cite{fang06}, we can
finish the proof of \textbf{Claim 1} in the following lemma.

\begin{lem}\label{ss-Lem3.4}
  Under the assumptions of Lemma \ref{ss-L2.2}, if $\epsilon_0$ is small enough, we obtain
        \begin{equation}
    |g(x,t)-g_\infty(x)|\leq C_{10}\epsilon_0^{\frac{\theta}{\gamma}2^{-N_4-1}},\label{ss-E2.21}
        \end{equation}
for all $x\in[0,M]$ and $t\in[0,T]$.
\end{lem}
\begin{proof}
From (\ref{ss-2.3}), for any fixed $x\in[0,M]$, we have
    \begin{eqnarray}
      &&I_1(x,t)+\frac{\theta}{2c_1+c_2}\int^t_0
      Ar^\beta(x,\tau)(\rho^\gamma(x,\tau)-\rho_\infty^\gamma(x))d\tau\nonumber\\
      &=&r_0^\beta(x)\rho_0^\theta(x)+
      I_2(x,t),\ x\in[0,M],\ t\in[0,T],
    \end{eqnarray}
where
    \begin{eqnarray*}
    &&I_1(x,t)\\
    &=&r_\infty^\beta(x)
      \rho^\theta(x,t)-(r_\infty^\beta(x)-r^\beta(x,t))\rho^\theta(x,t)
    +\int^M_x\beta[(r^{\beta-n}\rho^{\theta-1})(y,t)
    -(r_0^{\beta-n}\rho_0^{\theta-1})(y)]dy\nonumber\\
    &&-\frac{\theta}{2c_1+c_2}\int^M_x
    [(r^{\beta-n+1}u)(y,t)
    -(r_0^{\beta-n+1}u_0)(y)]dy
    +\frac{\theta(\beta-n+1)}{2c_1+c_2}\int^t_0\int^M_x
    r^{\beta-n}u^2dyd\tau,
    \end{eqnarray*}
and
    \begin{eqnarray*}
      &&I_2(x,t)\\
        &=&-\frac{\theta
        A\beta}{2c_1+c_2}\int^t_0\int^M_xr^{\beta-n}\frac{\rho^\gamma-\rho_\infty^\gamma}{\rho}dyd\tau\\
        &&+\frac{\theta}{2c_1+c_2}\int^t_0\int^M_x
        \left\{r^\beta
        G(M_0+y)(r^{2-2n}-r_\infty^{2-2n})+r^{\beta-n+1}\Delta f
        \right\}dyd\tau.
    \end{eqnarray*}

Using (\ref{ss-E2.9})-(\ref{ss-E2.10}), (\ref{ss-E3.66}),
H\"{o}lder's inequality and the condition
$\epsilon_2<\frac{\gamma-\theta-1}{\gamma-\theta}$, i.e.,
$\frac{\theta+1+(\gamma-\theta)\epsilon_2}{\gamma}<1$, we have
    \begin{eqnarray}
      |(r-r_\infty)(x)|&\leq& C|r^n-r_\infty^n|\leq
      C\int^x_0|\rho^{-1}-\rho_\infty^{-1}|dy\leq
            C\int^x_0(M-y)^{-\frac{1}{\gamma}}|g-g_\infty|dy\nonumber\\
      &\leq&
      C\left(\int^x_0(M-y)^\frac{\theta-1+(\gamma-\theta)\epsilon_2}{\gamma}(g-g_\infty)^2dy
      \right)^{\frac{1}{2}}\left(\int^x_0(M-y)^{-\frac{\theta+1+(\gamma-\theta)\epsilon_2}{\gamma}}dy
      \right)^\frac{1}{2}\nonumber\\
      &\leq&
      C\epsilon_0^{2^{-N_4-1}},\label{ss-E3.69}
    \end{eqnarray}
            \begin{eqnarray}
             && \int^M_x|\rho^{\theta-1}-\rho_\infty^{\theta-1}|dy
              \leq
            C\int^M_x(M-y)^{\frac{\theta-1}{\gamma}}|g-g_\infty|dy\nonumber\\
      &\leq&
      C\left(\int^M_x(M-y)^\frac{\theta-1+(\gamma-\theta)\epsilon_2}{\gamma}(g-g_\infty)^2dy
      \right)^{\frac{1}{2}}\left(\int^M_x(M-y)^{\frac{2\theta}{\gamma}-\frac{\theta+1+(\gamma-\theta)\epsilon_2}{\gamma}}dy
      \right)^\frac{1}{2}\nonumber\\
      &\leq&
      C\epsilon_0^{2^{-N_4-1}}(M-x)^{\frac{\theta}{\gamma}},
            \end{eqnarray}
using the fact $\theta\in(0,\frac{\gamma}{2}]$ and the estimate
(\ref{ss-E2.10})-(\ref{ss-E2.11}), we have
    \begin{eqnarray}
      &&\left|-\frac{\theta}{2c_1+c_2}\int^M_x
    [(r^{\beta-n+1}u)(y,t)
    -(r_0^{\beta-n+1}u_0)(y)]dy\right|\nonumber\\
    &\leq&C(M-x)^\frac{1}{2}(\|u\|_{L^2_x}+\|u_0\|_{L^2_x})\nonumber\\
    &\leq&C\epsilon_0^{2^{-N_4-1}}(M-x)^\frac{\theta}{\gamma},\ x\in[0,M].\label{ss-E3.70}
    \end{eqnarray}

 Thus, from (\ref{ss-epsilon1}),
(\ref{ss-E2.9})-(\ref{ss-E2.13}) and
(\ref{ss-E3.69})-(\ref{ss-E3.70}), we obtain
    \begin{equation}
    |I_1(x,t)-r_\infty^\beta\rho^\theta|\leq C_{1,1}(M-x)^{\frac{\theta}{\gamma}}
    \epsilon_0^{2^{-N_4-1}},
    \label{ss-E3.24}
    \end{equation}
and
    \begin{equation}
      |I_2(x,t_1)-I_2(x,t_2)|\leq
      C_{1,2}(M-x)\epsilon_0^{2^{-N_4-1}}|t_2-t_1|,\ x\in[0,M].\label{ss-E3.25}
    \end{equation}
\textbf{Claim 1.1.} For any fixed $x\in[0,M]$, we have
    $$
    I_1(x,t)\geq\min\left\{
    I_1(x,0),r^{\beta}_\infty\left(\rho_\infty^\gamma-\frac{C_{1,2}}{C_{1,3}}\epsilon_0^{2^{-N_4-1}}
    (M-x)
    \right)^{
    \frac{\theta}{\gamma}}-C_{1,1}\epsilon_0^{2^{-N_4-1}}(M-x)^{\frac{\theta}{\gamma}}
    \right\}:=M_{1,1},
    $$
where $C_{1,3}:=\frac{\theta a^\beta}{2c_1+c_2}\leq \frac{\theta
r^\beta}{2c_1+c_2}$.

\noindent\textbf{Proof of Claim 1.1.} If not, there exists
$t_{1,1}>0$ such that $I_1(x,t_{1,1})<M_{1,1}$, then we can find
$t_{1,2}\in(0,t_{1,1})$ such that $I_1(x,t_{1,2})=M_{1,1}$ and
$I_1(x,t)<M_{1,1}$ for all $t\in(t_{1,2},t_{1,1})$. From
(\ref{ss-E3.25}) we have
    $$
    I_1(x,t_{1,1})-I_1(x,t_{1,2})+\frac{\theta}{2c_1+c_2}\int^{t_{1,1}}_{t_{1,2}}r^\beta
    (\rho^\gamma-\rho_\infty^\gamma)\geq-C_{1,2}\epsilon_0^{2^{-N_4-1}}(M-x)(t_{1,1}-t_{1,2}).
    $$

From (\ref{ss-E3.24}), we have
    \begin{eqnarray*}
   && \rho^\theta(x,t)=r_\infty^{-\beta}(I_1(x,t)-(I_1(x,t)-r_\infty^\beta\rho^\theta))\\
   &\leq& r_\infty^{-\beta}(M_{1,1}+C_{1,1}\epsilon_0^{2^{-N_4-1}}(M-x)^{\frac{\theta}{\gamma}})
   \leq
   \left(\rho_\infty^\gamma-\frac{C_{1,2}}{C_{1,3}}\epsilon_0^{2^{-N_4-1}}(M-x)
   \right)^{\frac{\theta}{\gamma}},
    \end{eqnarray*}
and
    $$
    \rho^\gamma\leq
    \rho_\infty^\gamma-\frac{C_{1,2}}{C_{1,3}}\epsilon_0^{2^{-N_4-1}}(M-x),
    $$
then $I_1(x,t_{1,1})\geq I_1(x,t_{1,2})$. It is a contradiction.
Thus, \textbf{Claim 1.1} holds.

Similarly, we can obtain the following Claim.

\noindent \textbf{Claim 1.2.} For any fixed $x\in[0,M]$, we have
    $$
    I_1(x,t)\leq\max\left\{
    I_1(x,0),r^{\beta}_\infty\left(\rho_\infty^\gamma+\frac{C_{1,2}}{C_{1,4}}
    \epsilon_0^{2^{-N_4-1}}(M-x)
    \right)^{
    \frac{\theta}{\gamma}}+C_{1,1}\epsilon_0^{2^{-N_4-1}}(M-x)^{\frac{\theta}{\gamma}}
    \right\}:=M_{1,2},
    $$
where $C_{1,4}$ is a positive constant satisfying $C_{1,4}\geq
\frac{\theta r^\beta}{2c_1+c_2}$.

From \textbf{Claim 1.1} and \textbf{1.2}, we have
    $$
    |g(x,t)-g_\infty(x)|\leq
    C_{1,5}\epsilon_0^{\frac{\theta}{\gamma}2^{-N_4-1}},
    $$
where $x\in[0,M]$ and $t\in[0,T]$, when $\epsilon_0\leq
\delta_{11}$ is small enough.
\end{proof}

Now, we can let
    \begin{equation}
      \epsilon_1=\epsilon_0+C_{10}\epsilon_0^{\frac{\theta}{\gamma}2^{-N_4-1}}.\label{ss-self-E3.76}
    \end{equation}
If $4\epsilon_1<\min_{x\in[0,M]}g_\infty$, $C_3\epsilon_1\leq
\delta_9$, $\epsilon_1\leq \delta_{10}$ and $\epsilon_0\leq
\delta_{11}$, using the results in Lemmas
\ref{ss-L2.2}-\ref{ss-Lem3.4}, we finish the proof of the
\textbf{Claim 1}. From (\ref{ss-E3.5}) and Claim 1,
 using the classical continuity method, we can obtain  the
following lemma easily.

\begin{lem}\label{ss-L3.5}
  Under the assumptions of Theorem \ref{ss-thm},  we obtain
  (\ref{ss-E2.9})-(\ref{ss-E2.10}),
   (\ref{ss-E4.1-0})-\ref{ss-E4.2-3}), (\ref{ss-E2.21}) and
   \begin{equation}
     |r(x,t)-r_\infty(x)|\leq
     C_{11}\epsilon_0^{\frac{\theta}{\gamma}2^{-N_4-1}}\label{ss-E3.56}
   \end{equation}
hold for all $x\in[0,M]$ and $t\in[0,T^*)$.
\end{lem}
\begin{proof}
Let $\mathcal{A}=\{T\in[0,T^*)\left|\ I(t)\leq \epsilon_1\
\textrm{ for all }\ t\in[0,T]\right.\}$. Since $I(0)\leq
\epsilon_0<\epsilon_1$ and $I(t)\in C([0,T^*))$, then there exists
a constant $T_0>0$ such that $I(t)\leq \epsilon_1$ for all
$t\in[0,T_0]$. Thus,  $\mathcal{A}$ is not empty and relatively
closed in $[0,T^*)$.  To show that $\mathcal{A}$ is also
relatively open in $[0,T^*)\cap[0,T]$, and hence the entire
interval,
    it therefore suffices to show that the weaker bound
    $$
    I(t)\leq 2\epsilon_1,\ \textrm{ for all }\ t\in
    [0,T']\subset[0,T^*),
    $$
implies $I(t)\leq \epsilon_1$ for all $t\in[0,T']$. From
\textbf{Claim 1}, we have $\mathcal{A}=[0,T^*)$.

Then, from Lemmas \ref{ss-L2.2}-\ref{ss-Lem3.4}, we obtain
(\ref{ss-E2.9})-(\ref{ss-E2.10}),
   (\ref{ss-E4.1-0})-\ref{ss-E4.2-3}), (\ref{ss-E2.21}) and
   (\ref{ss-E3.56})
hold for all $x\in[0,M]$ and $t\in[0,T^*)$.
\end{proof}

We will prove an estimate in weighted $L^2([0,M]\times[0,T^*))$
norm of the function $g-g_\infty$.

\begin{lem}\label{ss-L3.10}
   Under the assumptions of Theorem \ref{ss-thm}, we obtain
    \begin{equation}
        \int^{t}_0\int^M_0(1+s)^{-\epsilon_2}
    (g-g_\infty)^2dxds\leq
    C,\label{ss-E4.4-10}
        \end{equation}
   where $t\in[0,T^*)$.
\end{lem}
\begin{proof}
From (\ref{ss-E}), we have
    \begin{eqnarray*}
      A(\rho^\gamma-\rho_\infty^\gamma)&=&\int^M_x\left(\frac{u_t}{r^{n-1}}+\frac{\Delta f}{r^{n-1}}
      \right)dy+\int^M_xG(M_0+y)(r^{2-2n}-r_\infty^{2-2n})dy\\
        &&+(2c_1+c_2)\rho^{\theta+1}(r^{n-1}u)_x-2c_1(n-1)\rho^\theta\frac{u}{r}
        -2c_1(n-1)\int^M_x\rho^\theta\left(\frac{u}{r}
        \right)_xdy
    \end{eqnarray*}
  Multiplying the above equality by $(1+t)^{-\epsilon_2} (M-x)^{-2}(\rho^\gamma-\rho_\infty^\gamma)$, integrating the
resulting equation over $[0,M]\times[0,t]$, we obtain
   \begin{eqnarray}
      &&\int^{t}_0\int^M_0A(1+s)^{{-\epsilon_2}}(M-x)^{-2}(\rho^\gamma-\rho_\infty^\gamma)^2dxds\nonumber\\
            &=&\int^{t}_0\int^M_0(1+s)^{{-\epsilon_2}}(M-x)^{-2}(\rho^\gamma-\rho_\infty^\gamma)\int^M_x
            \frac{u_t}{r^{n-1}}
     dydxds\nonumber\\
            &&+\int^{t}_0\int^M_0(1+s)^{{-\epsilon_2}}(M-x)^{-2}(\rho^\gamma-\rho_\infty^\gamma)\int^M_x
            \frac{\Delta f}{r^{n-1}}
     dydxds\nonumber\\
            &&+\int^t_0\int^M_0(1+s)^{{-\epsilon_2}}(M-x)^{-2}(\rho^\gamma-\rho_\infty^\gamma)
            \int^M_xG(M_0+y)(r^{2-2n}-r_\infty^{2-2n})dydxds\nonumber\\
            &&+\int^{t}_0\int^M_0(2c_1+c_2)(1+s)^{{-\epsilon_2}}(M-x)^{-2}(\rho^\gamma-\rho_\infty^\gamma)
            \rho^{\theta+1}(r^{n-1}u)_xdxds\nonumber\\
            &&-\int^{t}_0\int^M_02c_1(n-1)(1+s)^{{-\epsilon_2}}(M-x)^{-2}(\rho^\gamma-\rho_\infty^\gamma)
            \rho^\theta\frac{u}{r}dxds\nonumber\\
            &&-\int^{t}_0\int^M_02c_1(n-1)
            (1+s)^{{-\epsilon_2}}(M-x)^{-2}(\rho^\gamma-\rho_\infty^\gamma)\int^M_x\rho^\theta\left(\frac{u}{r}
        \right)_xdydxds\nonumber\\&
            :=&\sum^6_{i=1}F_i.\label{ss-E2.22-1-7}
    \end{eqnarray}

Using (\ref{ss-E}), (\ref{ss-epsilon1}),
(\ref{ss-E2.9})-(\ref{ss-E2.10}),
(\ref{ss-E4.1-0})-(\ref{ss-E4.2-3}),
 (\ref{ss-E2.21}), Lemma \ref{ss-L3.5},
integration by parts and the Cauchy-Schwarz inequality, we can
estimate $F_i$ as follows.
\begin{eqnarray}
  F_1&=&\left.\left\{\int^M_0(1+s)^{{-\epsilon_2}}(M-x)^{-2}(\rho^\gamma-\rho_\infty^\gamma)\int^M_x
            \frac{u}{r^{n-1}}
     dydx\right\}\right|^t_0\nonumber\\
        &&+\int^t_0\int^M_0{{\epsilon_2}}(1+s)^{{{-\epsilon_2}}-1}(M-x)^{-2}(\rho^\gamma-\rho_\infty^\gamma)\int^M_x
            \frac{u}{r^{n-1}}
     dydxds\nonumber\\
        &&+\int^t_0\int^M_0\gamma(1+s)^{{-\epsilon_2}}(M-x)^{-2}\rho^{\gamma+1}(r^{n-1}u)_x\int^M_x
            \frac{u}{r^{n-1}}
     dydxds\nonumber\\
        &&+\int^t_0\int^M_0(n-1)(1+s)^{{-\epsilon_2}}(M-x)^{-2}(\rho^\gamma-\rho_\infty^\gamma)\int^M_x
            \frac{u^2}{r^{n}}
     dydxds\nonumber\\
        &\leq&C\|g-g_\infty\|_{L^\infty_{xt}}\left(\int^M_0u^2dx\right)^{\frac{1}{2}}
        \left(\int^M_0(M-x)^{-\frac{1}{2}}dx
        \right)^\frac{1}{2}+C\nonumber\\
            &&+C\left(\int^t_0\int^M_x\frac{\rho_\infty^{\gamma-1}}{(1+s)^{\epsilon_2}}(g-g_\infty)^2dxds
            \right)^\frac{1}{2}\left(\int^t_0\|u(\cdot,s)\|_{L^\infty}^2ds\int^M_0(M-x)^{-\frac{\gamma-1}{\gamma}}dx
            \right)^\frac{1}{2}\nonumber\\
        &&+C\left(\int^t_0\int^M_x\rho^{\theta+1}(r^{n-1}u)_x^2dxds
            \right)^\frac{1}{2}\left(\int^t_0\|u(\cdot,s)\|_{L^\infty}^2ds\int^M_0(M-x)^{-\frac{\theta-1}{\gamma}}dx
            \right)^\frac{1}{2}
        \nonumber\\
                &&+C\|g-g_\infty\|_{L^\infty_{xt}}\int^t_0\|u(\cdot,s)\|_{L^\infty}^2ds\nonumber\\
      &\leq& C,
\end{eqnarray}
    \begin{equation}
      F_2\leq C\left(\int^t_0(1+t)^{-1{-\epsilon_2}}
      \right)^\frac{1}{2}
      \left(\int^t_0(1+t)f_1^2dt
      \right)^\frac{1}{2}\leq C,
    \end{equation}
        \begin{eqnarray}
       |r(x,t)-r_\infty(x)|&\leq&
        C\int^x_0\rho_\infty^{-1}|g-g_\infty|dy\leq C\left(\int^x_0\rho_\infty^{\gamma-1}(g-g_\infty)^2dy\right)^{\frac{1}{2}}
        \left(\int^x_0\rho_\infty^{-\gamma-1}dy\right)^{\frac{1}{2}}\nonumber\\
        &\leq&
        C(M-x)^{-\frac{1}{2\gamma}}\left(
        \int^M_0\rho_\infty^{\gamma-1}(g-g_\infty)^2dx\right)^{\frac{1}{2}},\label{ss-E3.62}
        \end{eqnarray}
    \begin{eqnarray}
      F_3
      &\leq&C\int^t_0\int^M_0(1+s)^{-\epsilon_2}(M-x)^{-2}|\rho^\gamma-\rho_\infty^\gamma|
      \left(
        \int^M_0\rho_\infty^{\gamma-1}(g-g_\infty)^2dx\right)^{\frac{1}{2}}\int^M_x(M-y)^{-\frac{1}{2\gamma}}dydx
        \nonumber\\
      &\leq&
      \frac{A}{4}\int^{t}_0\int^M_0(1+s)^{{-\epsilon_2}}(M-x)^{-2}(\rho^\gamma-\rho_\infty^\gamma)^2dxds\nonumber\\
      &&+C\int^{t}_0\int^M_0(1+s)^{{-\epsilon_2}}\rho_\infty^{\gamma-1}(g-g_\infty)^2dxds
      \int^M_0(M-z)^{-\frac{1}{\gamma}}dz\nonumber\\
       &\leq&
      \frac{A}{4}\int^{t}_0\int^M_0(1+s)^{{-\epsilon_2}}(M-x)^{-2}(\rho^\gamma
      -\rho_\infty^\gamma)^2dxds+C,
    \end{eqnarray}
        \begin{eqnarray}
          F_4&=&-\frac{2c_1+c_2}{\theta}\int^t_0\int^M_0(1+s)^{-\epsilon_2}(M-x)^{-2}(\rho^\gamma-\rho_\infty^\gamma)
          (\rho^\theta)_tdxds\nonumber\\
                &=&-\frac{2c_1+c_2}{\theta}\left.\left\{
                   (1+s)^{-\epsilon_2} \int^M_0(M-x)^{-2}\left(\frac{\theta}{\gamma+\theta}
                   \rho^{\gamma+\theta}-\rho_\infty^\gamma\rho^\theta\right)dx
                \right\}\right|^t_0\nonumber\\
          &&-\frac{\epsilon_2(2c_1+c_2)}{\theta}\int^t_0
                   \int^M_0(1+s)^{-\epsilon_2-1}(M-x)^{-2}\left(\frac{\theta}{\gamma+\theta}
                   \rho^{\gamma+\theta}-\rho_\infty^\gamma\rho^\theta\right) dx          \nonumber\\
          &\leq&C\|g-g_\infty\|_{L^\infty_{t,x}}\int^M_0(M-x)^{\frac{\theta}{\gamma}-1}dx
          \left(1+\int^t_0(1+s)^{-1-\epsilon_2}ds
          \right) +C\nonumber\\
                &\leq&C,
        \end{eqnarray}
    \begin{eqnarray}
      F_5&\leq&
      C\int^t_0\int^M_0(1+s)^{{-\epsilon_2}}\|u(\cdot,t)\|_{L^\infty}(M-x)^{\frac{\theta}{\gamma}-1}dxds\nonumber\\
            &\leq&
            C\left\{\int^t_0(1+s)^{-1-\epsilon_2}ds\right\}^\frac{1}{2}
            \left\{\int^t_0(1+s)^{1-\epsilon_2}\|u(\cdot,t)\|_{L^\infty}^2ds\right\}^\frac{1}{2}
            \leq             C
    \end{eqnarray}
and
   \begin{eqnarray}
      F_6 &\leq&
    C\int^t_0\int^M_0(1+s)^{{-\epsilon_2}}(M-x)^{-1}\int^M_x\left(|\rho^\theta
    u_x|+|\rho^{\theta-1}u|
    \right)dydxds\nonumber\\
        &\leq&C\int^t_0\int^M_0(1+s)^{{-\epsilon_2}}(M-x)^{-1}\left[\int^M_x\left(\rho^{\theta+1}
    u_x^2+\rho^{\theta-1}u^2
    \right)dy\right]^\frac{1}{2}\left[\int^M_x\rho^{\theta-1}dy\right]^\frac{1}{2}dxds\nonumber\\
        &\leq& C\int^t_0\int^M_0(1+s)^{{-\epsilon_2}}(M-x)^{\frac{\theta-1}{2\gamma}-\frac{1}{2}}\left[\int^M_x\left(\rho^{\theta+1}
    u_x^2+\rho^{\theta-1}u^2
    \right)dy\right]^\frac{1}{2}dxds\nonumber\\
        &\leq& C\left[\int^t_0\int^M_0(1+s)^{1-\epsilon_2}\left(\rho^{\theta+1}
    u_x^2+\rho^{\theta-1}u^2
    \right)dxds\right]^\frac{1}{2}
    \left[\int^t_0(1+s)^{-1-\epsilon_2}dt\right]^\frac{1}{2}\nonumber\\
        &\leq& C.\label{ss-E2.30-1-8}
    \end{eqnarray}
From (\ref{ss-E2.22-1-7})-(\ref{ss-E2.30-1-8}), we get
(\ref{ss-E4.4-10}) immediately.
\end{proof}
\begin{lem}\label{ss-L2.3}
 Under the assumptions of Theorem \ref{ss-thm}, we obtain
\begin{eqnarray}
          &&(1+t)^{-\epsilon_2}\int^M_0(M-x)^{1-\frac{\theta}{\gamma}}(\rho^\theta-\rho_\infty^\theta)_x^2(x,t)dx\nonumber\\
          &&+\int^t_0\int^M_0
          (1+s)^{-\epsilon_2}(M-x)^{2-\frac{2\theta}{\gamma}}
          (\rho^\theta-\rho_\infty^\theta)_x^2(x,s)dxds\leq C,
          \label{ss-E2.14-9}
        \end{eqnarray}
    for all $t\in[0,T^*)$.
\end{lem}
\begin{proof}From (\ref{ss-2.3}), we have
    \begin{eqnarray}
      &&\partial_t\left[\frac{\theta}{2c_1+c_2}r^{1+\beta-n}u
      +r^{\beta}(\rho^\theta)_x-r_\infty^{\beta}(\rho_\infty^\theta)_x
      \right]\nonumber\\
        &&+\frac{A\gamma\rho^{\gamma-\theta}}{2c_1+c_2}\left[\frac{\theta}{2c_1+c_2}r^{1+\beta-n}u
      +r^{\beta}(\rho^\theta)_x-r_\infty^{\beta}(\rho_\infty^\theta)_x
      \right]\nonumber\\
            &=&\frac{A\gamma\theta}{(2c_1+c_2)^2}\rho^{\gamma-\theta}r^{1+\beta-n}u
            +\frac{\theta(1+\beta-n)}{2c_1+c_2}r^{\beta-n}u^2
            -\frac{\theta}{2c_1+c_2}r^{1+\beta-n}\Delta
            f\nonumber\\
                &&-\frac{\theta}{2c_1+c_2}\left(
                r^\beta\frac{G(M_0+x)}{r^{2n-2}}+r_\infty^\beta\frac{\rho^{\gamma-\theta}}{\rho_\infty^{\gamma-\theta}}
                (A\rho_\infty^\gamma)_x                \right).
            \label{ss-E2.15}
    \end{eqnarray}
Let $H=\frac{\theta}{2c_1+c_2}r^{1+\beta-n}u
      +r^{\beta}(\rho^\theta)_x-r_\infty^{\beta}(\rho_\infty^\theta)_x
$, multiplying (\ref{ss-E2.15}) by
$(1+t)^{-\epsilon_2}(M-x)^{1-\frac{\theta}{\gamma}}H$, integrating
the resulting equation over $[0,M]$, and using the Cauchy-Schwarz
inequality, we obtain
    \begin{eqnarray}
      &&\frac{d}{dt}\int^M_0(1+t)^{-\epsilon_2}(M-x)^{1-\frac{\theta}{\gamma}}H^2(x,t)dx+C_{13}\int^M_0
      (1+t)^{-\epsilon_2}(M-x)^{2-\frac{2\theta}{\gamma}}H^2(x,t)dx\nonumber\\
      &\leq&C\int^M_0(1+t)^{-\epsilon_2}(M-x)^{1-\frac{\theta}{\gamma}}\left((M-x)^{1-\frac{\theta}{\gamma}}
      |Hu|+|u^2H|
      +|\Delta fH|\right)dx\nonumber\\
                &&+\int^M_0(1+t)^{-\epsilon_2}(M-x)^{1-\frac{\theta}{\gamma}}
                \left|r^\beta G\frac{M_0+x}{r^{2n-2}}+\frac{ \rho^{\gamma-\theta}r_\infty^{\beta}}{ \rho_\infty^{\gamma-\theta}}
                (A\rho_\infty^\gamma)_x\right||H|
                dx\nonumber\\
            &\leq& \frac{C_{13}}{4}\int^M_0
      (1+t)^{-\epsilon_2}(M-x)^{2-\frac{2\theta}{\gamma}}H^2(x,t)dx+C\|u(\cdot,t)\|_{L^\infty}^2+
      C\|u(\cdot,t)\|_{L^\infty}^2\|u(\cdot,t)\|_{L^2}^2\nonumber\\
              &&  +Cf_1^2
                  +C\int^M_0(1+t)^{-\epsilon_2}\left|r^\beta G\frac{M_0+x}{r^{2n-2}}+
                  \frac{ \rho^{\gamma-\theta}r_\infty^{\beta}}{ \rho_\infty^{\gamma-\theta}}
                (A\rho_\infty^\gamma)_x\right|^2dx.
            \label{ss-E2.16-9}
    \end{eqnarray}
Here, we use the estimates (\ref{ss-E2.9})-(\ref{ss-E2.10}) and
the condition $\theta\in(0,\gamma-1)$. From (\ref{ss-rhoinf3}) and
(\ref{ss-E2.9})-(\ref{ss-E2.10}),
 we have
    \begin{eqnarray}
   && \int^M_0(1+t)^{-\epsilon_2}\left|r^\beta G\frac{M_0+x}{r^{2n-2}}+
                  \frac{ \rho^{\gamma-\theta}r_\infty^{\beta}}{ \rho_\infty^{\gamma-\theta}}
                (A\rho_\infty^\gamma)_x\right|^2dx\nonumber\\
                    &=&\int^M_0(1+t)^{-\epsilon_2}\left|r^\beta
                    G\frac{M_0+x}{r^{2n-2}}-
                  G\frac{ (\rho^{\gamma-\theta}r_\infty^{\beta})(M_0+x)}{
                  ( \rho_\infty^{\gamma-\theta})(r_\infty^{2n-2})}\right|^2dx\nonumber\\
        &\leq&
        C\int^M_0(1+t)^{-\epsilon_2}\left[(r-r_\infty)^2+(g-g_\infty)^2
           \right]dx.\label{ss-E2.24-9}
    \end{eqnarray}
From (\ref{ss-E4.1-0}) and (\ref{ss-E2.16-9})-(\ref{ss-E2.24-9}),
we obtain
    \begin{eqnarray}
      &&\frac{d}{dt}\int^M_0(1+t)^{-\epsilon_2}(M-x)^{1-\frac{\theta}{\gamma}}
      H^2(x,t)dx+\frac{C_{13}}{2}\int^M_0(1+t)^{-\epsilon_2}(M-x)^{2-\frac{2\theta}{\gamma}}H^2(x,t)dx\nonumber\\
      &\leq& C\int^M_0(1+t)^{-\epsilon_2}\left(
            (r-r_\infty)^2+(g-g_\infty)^2
            \right)dx+C\|u(\cdot,t)\|_{L^\infty}^2+Cf_1^2.
            \label{ss-E2.36-9}
    \end{eqnarray}
From (A2), (\ref{ss-epsilon1}),
(\ref{ss-E4.1-0})-(\ref{ss-E4.2-3}),
(\ref{ss-E3.56})-(\ref{ss-E4.4-10}) and (\ref{ss-E2.36-9}), we
obtain (\ref{ss-E2.14-9}) immediately.

\end{proof}

\begin{lem}\label{ss-L3.11}
 Under the assumptions of Theorem \ref{ss-thm},  we obtain
            \begin{equation}
             (1+t)^{1-\epsilon_2} \int^M_0\left(\rho^{\theta-1}
              u^2+\rho^{\theta+1}u_x^2
              \right)(x,t)dx+\int^t_0\int^M_0(1+s)^{1-\epsilon_2}u_t^2(x,s)dxds\leq
              C,\label{ss-E2.17-9}
            \end{equation}
                \begin{equation}
                  \|u(\cdot,t)\|_{L^\infty}\leq C
                  (1+t)^{-\frac{1}{2}+\frac{\epsilon_2}{2}},\label{ss-E2.17-9-2}
                \end{equation}
  for all $t\in[0,T^*)$.
\end{lem}
\begin{proof}
   Multiplying (\ref{ss-E})$_2$ by $(1+t)^{1-\epsilon_2}u_t$, integrating the
resulting equation over $[0,M]\times[0,t]$, using integration by
parts and the boundary conditions (\ref{ss-Efd}), we obtain
    \begin{eqnarray}
      &&\int^t_0\int^M_0(1+s)^{1-\epsilon_2}u_t^2(x,s)dxds\nonumber\\
            &=&\int^t_0\int^M_0A(1+s)^{1-\epsilon_2}\rho^\gamma\partial_x(r^{n-1}u_t)dxds\nonumber\\
                &&-\int^t_0\int^M_0
      (2c_1+c_2)(1+s)^{1-\epsilon_2}\rho^{1+\theta}\partial_x(r^{n-1}u)
      \partial_x(r^{n-1}u_t)dxds\nonumber\\
      &&+\int^t_0\int^M_02c_1(n-1)(1+s)^{1-\epsilon_2}\rho^\theta\partial_x(r^{n-2}uu_t)
      dxds-\int^t_0\int^M_0(1+s)^{1-\epsilon_2}fu_tdxds\nonumber\\
            &:=&\sum^4_{i=1}H_i.\label{ss-E2.18-9}
    \end{eqnarray}
Using (A3), (\ref{ss-epsilon1}), (\ref{ss-E2.9})-(\ref{ss-E2.10}),
 (\ref{ss-E4.1-0})-(\ref{ss-E4.2-0}) and the Cauchy-Schwarz
inequality, we obtain
    \begin{eqnarray}
      &&H_2+H_3\nonumber\\
            &=&\left.\left\{(1+s)^{1-\epsilon_2}\int^M_0\left[-\frac{2c_1+c_2}{2}\rho^{1+\theta}[\partial_x(r^{n-1}u)]^2
            +c_1(n-1)\rho^\theta\partial_x(r^{n-2}u^2)\right]dx\right\}\right|^t_0\nonumber\\
                &&+\int^t_0\int^M_0
      \frac{(2c_1+c_2)(1-\epsilon_2)}{2}(1+s)^{-\epsilon_2}\rho^{1+\theta}(r^{n-1}u)_x^2
      dxds\nonumber\\
            &&-\int^t_0\int^M_0c_1(n-1)(1-\epsilon_2)(1+s)^{-\epsilon_2}\rho^\theta\partial_x(r^{n-2}u^2)
      dxds\nonumber\\
            &&+\int^t_0\int^M_0(1+s)^{1-\epsilon_2}\left\{
                    (2c_1+c_2)(n-1)\rho^{1+\theta}\partial_x(r^{n-1}u)
                    \partial_x(r^{n-2}u^2)\right.\nonumber\\
                    &&-
      \frac{(2c_1+c_2)}{2}(1+\theta)\rho^{2+\theta}[\partial_x(r^{n-1}u)]^3
      +2\theta
                c_1(n-1)\rho^{\theta+1}\frac{u}{r}[\partial_x(r^{n-1}u)]^2
                \nonumber\\
    &&-\theta
    c_1n(n-1)\rho^\theta\frac{u^2}{r^2}\partial_x(r^{n-1}u)+
            2nc_1(n-1)(n-2)\rho^{\theta-1}\frac{u^3}{r^3}\nonumber\\
    &&-\left.3c_1(n-1)(n-2)\rho^\theta
            \frac{u^2}{r^2}\partial_x(r^{n-1}u)\right\}dxds\nonumber\\
                        &\leq&C+C\int^t_0\left(\|u\|_{L^\infty_x}+\|\rho(r^{n-1}u)_x\|_{L^\infty_x}
    \right)(1+s)^{1-\epsilon_2}\int^M_0\left[\rho^{1+\theta}(r^{n-1}u)_x^2+\rho^{\theta-1}u^2\right]dxds\nonumber\\
    &&-C_{14}(1+t)^{1-\epsilon_2}
                        \int^M_0\left[\rho^{1+\theta}(r^{n-1}u)_x^2+\rho^{\theta-1}u^2\right]dx,
    \end{eqnarray}
        \begin{eqnarray}
          H_1
          &=&\left.\left\{(1+s)^{1-\epsilon_2}\int^M_0A\rho^\gamma\partial_x(r^{n-1}u)dx\right\}\right|^t_0
          +\int^t_0\int^M_0A\gamma(1+s)^{1-\epsilon_2}\rho^{\gamma+1}[\partial_x(r^{n-1}u)]^2dxds\nonumber\\
                        &&-\int^t_0\int^M_02A(n-1)(1+s)^{1-\epsilon_2}\rho^{\gamma}\frac{u}{r}\partial_x(r^{n-1}u)dxds
                        \nonumber\\
          &&
          +\int^t_0\int^M_0An(n-1)(1+s)^{1-\epsilon_2}\rho^{\gamma-1}\frac{u^2}{r^2}dxds \nonumber\\
          && -\int^t_0\int^M_0A(1-\epsilon_2)(1+s)^{-\epsilon_2}\rho^\gamma\partial_x(r^{n-1}u)dxds\nonumber\\
          &\leq&(1+s)^{1-\epsilon_2}\int^M_0A\rho^\gamma\partial_x(r^{n-1}u)dx+C,
        \end{eqnarray}
        \begin{eqnarray}
          H_4&=&-\left.\left\{(1+s)^{1-\epsilon_2}\int^M_0G\frac{u(M_0+x)}{r^{n-1}}dx
          \right\}\right|^t_0
          +(1-\epsilon_2)\int^t_0\int^M_0(1+s)^{-\epsilon_2}\frac{Gu(M_0+x)}{r^{n-1}}dxds\nonumber\\
            &&
          +\int^t_0\int^M_0(1-n)(1+s)^{1-\epsilon_2}G(M_0+x)r^{-n}u^2dxds-
          \int^t_0\int^M_0(1+s)^{1-\epsilon_2}\Delta fu_tdxds\nonumber\\
            &\leq& -(1+s)^{1-\epsilon_2}\int^M_0G
            \frac{u(M_0+x)}{r^{n-1}}dx+C+\frac{1}{2}\int^t_0\int^M_0(1+s)^{1-\epsilon_2}u_t^2dxds.
        \end{eqnarray}
Using (\ref{ss-E2.9})-(\ref{ss-E2.10}), (\ref{ss-E3.32}),
(\ref{ss-E4.1-0}), integration by parts and the Cauchy-Schwarz
inequality, we obtain
        \begin{eqnarray}
        &&(1+s)^{1-\epsilon_2}\int^M_0A\rho^\gamma\partial_x(r^{n-1}u)dx-(1+s)^{1-\epsilon_2}\int^M_0G
            \frac{u(M_0+x)}{r^{n-1}}dx\nonumber\\
                &=&(1+s)^{1-\epsilon_2}\int^M_0\left(A\rho^\gamma\partial_x(r^{n-1}u)+
            r^{n-1}u(A\rho_\infty^\gamma)_x-G
            r^{n-1}u(M_0+x)(r^{2-2n}-r_\infty^{2-2n})\right)dx\nonumber\\
                            &=&(1+s)^{1-\epsilon_2}\int^M_0\left(A(\rho^\gamma
                            -\rho_\infty^\gamma)\partial_x(r^{n-1}u)-G
                            r^{n-1}u(M_0+x)(r^{2-2n}-r_\infty^{2-2n})\right)dx\nonumber\\
                    &\leq&\frac{C_{14}}{4}(1+t)^{1-\epsilon_2}
                        \int^M_0\rho^{1+\theta}(r^{n-1}u)_x^2dx+C\int^M_0(1+t)^{1-\epsilon_2}
                        [\rho^{\gamma-1}_\infty(g-g_\infty)^2+(r-r_\infty)^2]dx\nonumber\\
                    &\leq&\frac{C_{14}}{4}(1+t)^{1-\epsilon_2}
                        \int^M_0\rho^{1+\theta}(r^{n-1}u)_x^2dx+C.\label{ss-E2.22-9}
        \end{eqnarray}
 From (\ref{ss-E2.18-9})-(\ref{ss-E2.22-9}),   we can obtain
\begin{eqnarray}
            &&(1+t)^{1-\epsilon_2}\int^M_0\left[\rho^{1+\theta}(r^{n-1}u)_x^2+\rho^{\theta-1}u^2\right]dx+
            \int^t_0\int^M_0(1+s)^{1-\epsilon_2}u_t^2dxds\nonumber\\
    &\leq& C+C  \int^t_0     \left(\|u\|_{L^\infty_x}+
    \|\rho(r^{n-1}u)_x\|_{L^\infty_{x}}\right)
           (1+s)^{1-\epsilon_2}\nonumber\\
                &&\times \int^M_0\left[\rho^{1+\theta}(r^{n-1}u)_x^2+\rho^{\theta-1}u^2\right]dx
           ds.
    \label{ss-E2.45}
    \end{eqnarray}
From (\ref{ss-2.4}), we have
    \begin{eqnarray*}
            &&\rho(r^{n-1}u)_x\nonumber\\
                    &=&\frac{1}{(2c_1+c_2)\rho^{\theta}}\left\{A\rho^\gamma+2c_1(n-1)\rho^\theta\frac{u}{r}
            +\int^M_x\left[-\frac{u_t}{r^{n-1}}
                +2c_1(n-1)\rho^\theta\left(\frac{u}{r}\right)_x-\frac{f}{r^{n-1}}\right]dy\right\}.
            \end{eqnarray*}
Using conditions $\theta\in(0,\frac{\gamma}{2}]\cap(0,\gamma-1)$,
(\ref{ss-epsilon1}), estimates (\ref{ss-E2.9})-(\ref{ss-E2.13})
and H\"{o}lder's inequality, we conclude that
    \begin{equation}
    \|\rho\partial_x(r^{n-1}u)\|_{L^\infty_{x}}\leq
              C+C \left(\|u(\cdot,t)\|_{L^\infty}^2+
             \int^M_0\left[\rho^{1+\theta}(r^{n-1}u)_x^2+u_t^2\right]dx\right)^\frac{1}{2}.
              \label{ss-E2.70}
    \end{equation}
Using (\ref{ss-E3.32}), (\ref{ss-E4.2-0}),
(\ref{ss-E2.45})-(\ref{ss-E2.70}) and the Cauchy-Schwarz
inequality, we can obtain
    \begin{eqnarray*}
            &&(1+t)^{1-\epsilon_2}\int^M_0\left[\rho^{1+\theta}(r^{n-1}u)_x^2+\rho^{\theta-1}u^2\right]dx+
            \int^t_0\int^M_0(1+s)^{1-\epsilon_2}u_t^2dxds\nonumber\\
    &\leq& C+C  \int^t_0      \left(
            \int^M_0(1+s)^{1-\epsilon_2}\left[\rho^{1+\theta}(r^{n-1}u)_x^2+\rho^{\theta-1}u^2\right]dx\right)^2
            ds.
       \end{eqnarray*}
Using Gronwall's inequality
 and the estimate (\ref{ss-E4.2-0}), we obtain
(\ref{ss-E2.17-9}) immediately.

From (\ref{ss-E2.9}), (\ref{ss-E2.17-9}) and   the fact
$\theta\in(0,\gamma-1)$, we can obtain
    \begin{eqnarray*}
    |u(x,t)|&\leq& \left|\int^x_0 u_xdy
    \right|\leq C\left(\int^x_0 \rho^{\theta+1}u_x^2dy
    \right)^\frac{1}{2}\left(\int^x_0 (M-y)^{-\frac{\theta+1}{\gamma}}dy
    \right)^\frac{1}{2}\nonumber\\
            &\leq& C
             (1+t)^{-\frac{1}{2}+\frac{\epsilon_2}{2}}, \
             (x,t)\in[0,M]\times[0,T^*).
    \end{eqnarray*}
\end{proof}

\begin{lem}\label{ss-L2.8}
  Under the assumptions of Theorem \ref{ss-thm}, we obtain
        \begin{equation}
          \int^M_0 u_t^2(x,t)dx+\int^t_0 \int^M_0 \left[
            \rho^{\theta+1}u_{xt}^2+\rho^{\theta-1}u_t^2
            \right]dxds\leq C_{11},\label{ss-E2.36}
        \end{equation}
            \begin{equation}
              \|\rho(r^{n-1}u)_x(\cdot,t)\|_{L^\infty}\leq
              C_{11},\label{ss-E2.37}
            \end{equation}
for all $t\in[0,T^*)$.
\end{lem}
\begin{proof}
  We differentiate the equation (\ref{ss-E})$_2$ with respect to
  $t$, multiply it by $u_t$ and integrate it over
  $[0,M]\times[0,t]$, using the boundary conditions
(\ref{ss-Efd}), then derive
    \begin{eqnarray}
      &&\int^M_0\frac{1}{2} u_t^2dx\nonumber\\
               &=&\int^M_0\frac{1}{2} u_t^2(x,0)dx-\int^t_0\int^M_0\left[(2c_1+c_2)\rho^{1+\theta}\partial_x(r^{n-1}u)
                -A\rho^\gamma
               -2c_1(n-1)\rho^\theta\frac{u}{r}\right]\nonumber\\
               &&\times\partial_x((n-1)r^{n-2}uu_t)dxds
    -\int^t_0\int^M_0\partial_t\left[(2c_1+c_2)\rho^{1+\theta}\partial_x(r^{n-1}u)
    -A\rho^\gamma\right.\nonumber\\
    &&\left.    -2c_1(n-1)\rho^\theta\frac{u}{r}\right]
    \partial_x(r^{n-1}u_t)dxds+\int^t_0\int^M_02c_1(n-1)\partial_t(r^{n-1}\rho^\theta\partial_x(\frac{u}{r}))
             u_tdxds\nonumber\\
            &&-\int^t_0\int^M_0f_t u_tdxds\nonumber\\
    &:=&\sum_{i=1}^5J_i.\label{ss-E2.38}
    \end{eqnarray}
    From (A2)-(A3), we have
        \begin{eqnarray}
          J_1&\leq& C\left(\left\|\left((2c_1+c_2)\rho_{0}^{\theta+1}(r_{0}^{n-1}u_{0})_{x}\right)_{x}
    -2c_1(n-1)\frac{u_0}{r_0}(\rho_0^{\theta})_x\right\|_{L^2}\right.
    \nonumber\\
        &&\left.+\|(\rho_0^\gamma)_x\|_{L^2}+\|f(x,r_0,0)\|_{L^2} \right)^2\nonumber\\
        &\leq& C.
        \end{eqnarray}
From (\ref{ss-E2.9})-(\ref{ss-E2.13}) and  the Cauchy-Schwarz
inequality, we get
        \begin{eqnarray}
                &&J_3+J_4\nonumber\\
          &=&-\int^t_0\int^M_0\left[(2c_1+c_2)\rho^{1+\theta}(r^{n-1}u_t)_x^2
          -2c_1(n-1)\rho^\theta(r^{n-2}u_t^2)_x
          \right]dxds\nonumber\\
                &&+\int^t_0\int^M_0\big\{(2c_1+c_2)(1+\theta)\rho^{\theta+2}
                [\partial_x(r^{n-1}u)]^2-(n-1)(2c_1+c_2)\rho^{1+\theta}\partial_x(r^{n-2}u^2)
                \nonumber\\
          &&-\gamma\rho^{\gamma+1}\partial_x(r^{n-1}u)-2c_1(n-1)\theta\rho^{\theta+1}
          \partial_x(r^{n-1}u)\frac{u}{r}-2c_1(n-1)\rho^\theta\frac{u^2}{r^2}
                \big\}\nonumber\\
                &&\times[(n-1)\frac{u_t}{r\rho}
                +r^{n-1}u_{tx}]dxds
                +2c_1(n-1)\int^t_0\int^M_0\big\{(n-1)r^{n-2}u\rho^\theta\left(\frac{u}{r}
                \right)_xu_t
  \nonumber\\
            &&-\theta
  r^{n-1}\rho^{\theta+1}(r^{n-1}u)_x\left(\frac{u}{r}
  \right)_xu_t
  -r^{n-1}\rho^\theta\left(\frac{u^2}{r^2}
  \right)_xu_t\big\}dxds
            \nonumber\\
  &\leq&
  -C_{15}\int^t_0\int^M_0(\rho^{\theta+1}u_{xt}^2+\rho^{\theta-1}u_t^2)dxds+C
\nonumber\\
            && +C\int^t_0 \left(\|u\|_{L^\infty_x}^2+
            \|\rho(r^{n-1}u)_x\|_{L^\infty_{x}}^2
            \right)\int^M_0 \left[
            \rho^{\theta+1}u_x^2+\rho^{\theta-1}u^2
            \right]dxds,
        \end{eqnarray}
     From  (\ref{ss-f4}), (\ref{ss-E2.9})-(\ref{ss-E2.13})
and the Cauchy-Schwarz inequality, we obtain
            \begin{eqnarray}
              J_2
               &\leq& \frac{C_{15}}{8} \int^t_0\int^M_0 \left(\rho^{\theta-1}u_t^2+
               \rho^{\theta+1}u_{xt}^2\right)dxds+C\nonumber\\
            && +C\int^t_0 \left(\|u\|_{L^\infty_x}^2+
            \|\rho(r^{n-1}u)_x\|_{L^\infty_{x}}^2
            \right)\int^M_0 \left[
            \rho^{\theta+1}u_x^2+\rho^{\theta-1}u^2
            \right]dxds
            \end{eqnarray}
and
             \begin{eqnarray}
      J_5&\leq
      &\frac{C_{15}}{8}\int^t_0\int^M_0\rho^{\theta-1}u_t^2dxds\nonumber\\
        &&      +C\int^t_0\int^M_0(G(M_0+x)r^{-n}|u|
      +|\partial_r\Delta fu|+|\partial_t\Delta f|)^2\rho^{1-\theta}dxds\nonumber\\
                &\leq &\frac{C_{15}}{8}\int^M_0r^{\alpha-2}u_t^2dxds
              +C.\label{ss-E2.40}
    \end{eqnarray}
 From
(\ref{ss-E2.38})-(\ref{ss-E2.40}), we have
         \begin{eqnarray}
          &&\int^M_0 u_t^2(x,t)dx+\int^t_0 \int^M_0 \left[
            \rho^{\theta+1}u_{xt}^2+\rho^{\theta-1}u_t^2
            \right]dxds\nonumber\\
                &\leq& C+C\int^t_0 \left(\|u\|_{L^\infty_x}^2+
            \|\rho(r^{n-1}u)_x\|_{L^\infty_{x}}^2
            \right)\int^M_0 \left[
            \rho^{\theta+1}u_x^2+\rho^{\theta-1}u^2
            \right]dxds.
          \label{ss-E2.68}
    \end{eqnarray}
From (\ref{ss-E2.17-9})-(\ref{ss-E2.17-9-2}) and (\ref{ss-E2.70}),
we have
    \begin{equation}
  \|\rho\partial_x(r^{n-1}u)\|_{L^\infty_{x}}\leq
              C+C\|u_t\|_{L^2_x}.
              \label{ss-E2.69}
    \end{equation}
From (\ref{ss-E2.13}),  (\ref{ss-E2.17-9})-(\ref{ss-E2.17-9-2})
and  (\ref{ss-E2.68})-(\ref{ss-E2.69}), we obtain
    \begin{eqnarray*}
          &&\int^M_0 u_t^2(x,t)dx+\int^t_0 \int^M_0 \left[
            \rho^{\theta+1}u_{xt}^2+\rho^{\theta-1}u_t^2
            \right]dxds\nonumber\\
                &\leq& C+C\int^t_0 \int^M_0 \left[
            \rho^{\theta+1}u_x^2+\rho^{\theta-1}u^2
            \right]dx\|u_t\|_{L^2_x}^2ds.
    \end{eqnarray*}
Using Gronwall's inequality and the estimate (\ref{ss-E2.13}), we
obtain (\ref{ss-E2.36})-(\ref{ss-E2.37}) immediately.
\end{proof}

\noindent \textbf{Proof of existence and uniqueness}

 If
$T^*<\infty$, from Lemmas \ref{ss-L3.5}-\ref{ss-L2.8}, we have,
for all $t\in[0,T^*)$,
    $$C^{-1}(M-x)^{1/\gamma}\leq\rho(x,t)\leq
    C(M-x)^{1/\gamma},\ x\in[0,M],$$
    $$\int^M_0(M-x)^{1-\frac{\theta}{\gamma}}(\rho^\theta)_x^2dx\leq C,
        \ \int^M_0(\rho^{\gamma})_x^2dx\leq C,$$
    $$\int^M_0u^2+(M-x)^{\frac{\theta+1}{\gamma}}u_x^2dx\leq C, \ u(0,t)=0,$$
    $$\int^M_0\left\{
    \left((2c_1+c_2)\rho^{\theta+1}(r^{n-1}u)_{x}\right)_{x}
    -2c_1(n-1)\frac{u}{r}\partial_{x}\rho^{\theta}\right\}^2dx\leq
    C.
    $$
Thus, from Theorem \ref{ss-localthm}, there exists $T_2>0$ such
that the free boundary problem (\ref{ss-E})-(\ref{ss-Efd}) admits
a unique weak solution $(\rho_2, u_2, r_2)(x,t)$ on $[0,
M]\times[T^*-\frac{T_2}{2}, T^*+\frac{T_{2}}{2}]$, with initial
data $(\rho,u,r)(x,T^*-\frac{T_2}{2})$. Using the uniqueness
result in Theorem \ref{ss-localthm}, we obtain that
    $$
    (\tilde{\rho},\tilde{u},\tilde{r})(x,t)=\left\{
    \begin{array}{ll}
    (\rho,u,r)(x,t), &t\in[0,T^*-\frac{T_2}{2}]\\
     (\rho_2, u_2, r_2)(x,t),  &t\in[T^*-\frac{T_2}{2}, T^*+\frac{T_2}{2}]
    \end{array}
    \right.
    $$
is a solution of the system (\ref{ss-E})-(\ref{ss-Efd}), which is
contradiction with the definition of $T^*$. Thus, we have
$T^*=\infty$. From Lemma \ref{ss-L3.5}-\ref{ss-L2.8}, we can show
that the global weak solution satisfies the regularity conditions
(\ref{ss-E1.32})-(\ref{ss-E1.32-7}) and (\ref{ss-E1.46}) in
Theorem \ref{ss-thm}.
\begin{rem}
  The uniqueness of the solution of Theorem
  \ref{ss-localthm} is obtained by the energy method.
 Let  $(u_i,\rho_i,r_i)$, $i=1,2$, be two solutions of the
    system (\ref{ss-E})-(\ref{ss-Efd}) satisfying the regularity conditions in Theorem \ref{ss-thm}.
    Using similar arguments as that in the uniqueness part in \cite{ChenP},  we can obtain, for all $T>0$,
        \begin{eqnarray*}
      &&\frac{d}{dt}\int^M_0
        \left(
         w^2
         +\rho_1^{1-\theta}\rho_2^{2\theta-4}\varrho^2
         +\rho_1^{\theta}\rho_2^{-1}\mathcal{R}^2
        \right)dx\\
      &&+C^{-1}\int^M_0\rho_1^{1+\theta}\left(\rho_1r^{2n-2}_1\left(\partial_xw\right)^2+
      \frac{w^2}{r^2_1\rho_1}\right)dx\\
        &\leq& C\int^M_0
        \left(
        w^2
        +\rho_1^{1-\theta}\rho_2^{2\theta-4}\varrho^2
        +\rho_1^{\theta}\rho_2^{-1}\mathcal{R}^2
        \right)dx, \ t\in[0,T],
    \end{eqnarray*}
where $(w,\varrho,\mathcal{R})=(u_1-u_2,\rho_1-\rho_2,r_1-r_2)$.
  Using Gronwall's inequality,  we could obtain
  $(u_1,\rho_1,r_1)=(u_2,\rho_2,r_2)$, a.e. $(x,t)\in[0,M]\times[0,T]$.
\end{rem}

\section{Further decay result}\label{ss-Sec4}

\begin{lem}\label{ss-L3.13}
 Let $\nu$ be a positive constant satisfying $\nu<\min\{1,\frac{2\gamma-2}{\gamma+\theta}\}$.
   Under the assumptions of Theorem \ref{ss-thm},  we obtain
            \begin{equation}
                \left\|\rho^\frac{\gamma+\theta}{2}(\cdot,t)-\rho_\infty^\frac{\gamma+\theta}{2}(\cdot)
                \right\|_{L^\infty}\leq
             C(1+t)^{-\frac{1}{4}+\frac{\epsilon_2}{2}},\label{ss-E3.92-1}
            \end{equation}
                \begin{equation}
                  \|r(\cdot,t)-r_\infty(\cdot)\|_{L^\infty}\leq
                  C_\nu
                  (1+t)^{-\frac{1}{4}\nu+\frac{\epsilon_2\nu}{2}}, \label{ss-E3.92-2}
                \end{equation}
 for all $t\geq0$, where $C_\nu$ is a positive constant depending on $\nu$.
\end{lem}
\begin{proof}
  From (\ref{ss-E2.21}), (\ref{ss-E3.56})  and (\ref{ss-E2.14-9}), we have
    \begin{equation}
      \int^M_0\left(\rho^\frac{\gamma+\theta}{2}-\rho_\infty^\frac{\gamma+\theta}{2}\right)_x^2dx\leq
      C(1+t)^{{\epsilon_2}}, \ t\geq0.
    \end{equation}
Combining (\ref{ss-E4.1-0}) and the Galiardo-Nirenberg inequality
$\|\phi\|_{L^\infty}\leq
\|\phi\|_{L^2}^\frac{1}{2}\|\phi'\|_{L^2}^\frac{1}{2}$, we obtain
    $$
    \left\|\rho^\frac{\gamma+\theta}{2}-\rho_\infty^\frac{\gamma+\theta}{2}
    \right\|_{L^\infty}\leq
             C(1+t)^{-\frac{1}{4}+\frac{\epsilon_2}{2}}, \ t\geq0.
    $$
From (\ref{ss-E2.9})-(\ref{ss-E2.10}), (\ref{ss-E2.21}) and
(\ref{ss-E3.92-1}), we have
    $$
     \left\|\rho^{\frac{\nu(\gamma+\theta)}{2}}(\cdot,t)-\rho_\infty^{\frac{\nu(\gamma+\theta)}{2}}(\cdot)
     \right\|_{L^\infty}
    \leq C\|\rho^\frac{\gamma+\theta}{2}(\cdot,t)-\rho_\infty^\frac{\gamma+\theta}{2}(\cdot)\|_{L^\infty}^{\nu}\leq
             C(1+t)^{-\frac{\nu}{4}+\frac{\epsilon_2\nu}{2}},
    $$
        \begin{eqnarray*}
     \|r(\cdot,t)-r_\infty(\cdot)\|_{L^\infty}&\leq&
     C\int^M_0|\rho^{-1}-\rho_\infty^{-1}|dx\nonumber\\
     &\leq&
     C\|\rho^{\frac{\nu(\gamma+\theta)}{2}}(\cdot,t)-\rho_\infty^{\frac{\nu(\gamma+\theta)}{2}}
     (\cdot)\|_{L^\infty}\int^M_0
     (M-x)^{-\frac{1}{\gamma}-\frac{\nu(\gamma+\theta)}{2\gamma}}dx\nonumber\\
            &\leq&C_\nu(1+t)^{-\frac{\nu}{4}+\frac{\epsilon_2\nu}{2}},
       \end{eqnarray*}
    for all $t\geq0$.
\end{proof}

Using similar arguments as that in Lemmas \ref{ss-L3.7},
\ref{ss-L2.3}-\ref{ss-L3.11} and \ref{ss-L3.13} with
$\nu=\frac{\gamma-1}{2\gamma}$, we can obtain the following lemma
and omit the proof.

\begin{lem}\label{ss-L3.14}
 Under the assumptions of Theorem \ref{ss-thm}, we have
    \begin{equation}
   \int^M_0(M-x)^\frac{\theta-1}{\gamma}(g-g_\infty)^2dx+ \int^{t}_0\int^M_0
   \left[ \rho_\infty^{\gamma-1}
    (g-g_\infty)^2+(r-r_\infty)^2\right]dxds\leq
    C,\label{ss-E3.92}
        \end{equation}
  \begin{equation}
      \int^M_0\left\{u^2(x,t)+(M-x)^{\frac{\gamma-1}{\gamma}}(g-g_\infty)^2
      +(r-r_\infty)^2\right\}dx\leq
      C(1+t)^{-1},\label{ss-E3.90}
        \end{equation}
        \begin{equation}
    \int^t_0(1+s)\|u(\cdot,s)\|_{L^\infty}^2ds+    \int^t_0\int^M_0
    (1+s)\left(\rho^{\theta+1}u_x^2+\rho^{\theta-1}u^2
        \right)(x,s)dxds\leq C,
        \end{equation}
        \begin{equation}
          \int^M_0(M-x)^{1-\frac{\theta}{\gamma}}(\rho^\theta-\rho_\infty^\theta)_x^2(x,t)dx+\int^t_0\int^M_0
          (M-x)^{2-\frac{2\theta}{\gamma}}
          (\rho^\theta-\rho_\infty^\theta)_x^2(x,s)dxds\leq C,\label{ss-E3.93}
        \end{equation}
            \begin{equation}
             (1+t) \int^M_0\left(\rho^{\theta-1}
              u^2+\rho^{\theta+1}u_x^2
              \right)(x,t)dx+\int^t_0\int^M_0(1+s)u_t^2(x,s)dxds\leq
              C,
            \end{equation}
 for all $t\geq0$.
\end{lem}
\begin{rem}
The key point is: similar to (\ref{ss-E2.30-1}), using the
estimates (\ref{ss-E2.10}), (\ref{ss-E4.2-0}), (\ref{ss-E3.92-2})
and the condition $\epsilon_2<\frac{\gamma-1}{2(3\gamma-1)}$, we
have
    \begin{eqnarray*}
      &&\int^t_0\int^M_0            2c_1(n-1)\rho^\theta \left(\frac{u}{r}(\frac{r^n}{n}-\frac{r^n_\infty}{n})
                   \right)_xdxds\\
                    &\leq&C\int^t_0\int^M_0\left\{|r-r_\infty|\left(|\rho^{\theta}u_x|+|\rho^{\theta-1}u|
                    \right)+\rho^\theta |u(\rho^{-1}-\rho_\infty^{-1})|\right\}dxds\\
      &\leq& C\int^t_0\int^M_0
    (1+s)^{1-\epsilon_2}\left(\rho^{\theta+1}u_x^2+\rho^{\theta-1}u^2
        \right)(x,s)dxds\\
                &&+C\int^t_0
        (1+s)^{\epsilon_2-1}\|r(\cdot,t)-r_\infty(\cdot)\|_{L^\infty}^2\int^M_0(M-x)^{\frac{\theta-1}{\gamma}}dxds\\
            &&+C\int^t_0
        (1+s)^{\epsilon_2-1}\|\rho^{\frac{\nu(\gamma+\theta)}{2}}(\cdot,t)-\rho_\infty^{\frac{\nu(\gamma+\theta)}{2}}
     (\cdot)\|_{L^\infty}^2\int^M_0(M-x)^{\frac{\theta+1}{\gamma}-\frac{2}{\gamma}-\frac{\nu(\gamma+\theta)}{\gamma}}dxds\\
            &\leq&C+C\int^t_0(1+s)^{\epsilon_2-1-\frac{\nu}{2}+{\epsilon_2\nu}}ds\leq C,\ t\geq0.
    \end{eqnarray*}
\end{rem}

Without loss of generality, we assume
$\eta\in(0,\frac{2(\gamma+\theta)}{\gamma+\theta+1})$.  Let
$\epsilon_4\in(0,\frac{\gamma+\theta-1}{3(\gamma+\theta)})$ be a
constant satisfying
$\frac{1-\epsilon_4}{2-\frac{3\gamma+3\theta-1}{2(\gamma+\theta)}+\frac{\epsilon_4}{2}}>\frac{2(\gamma+\theta)}{\gamma+\theta+1}-\eta$.
Define $\{\kappa_j\}$ and $\{\eta_j\}$ by $\eta_{j+1}=1+\kappa_j$,
$\kappa_j=\frac{3\gamma+3\theta-1}{4(\gamma+\theta)}\eta_j-\frac{\epsilon_4}{4}\eta_j-\frac{1}{2}-\frac{\epsilon_4}{2}$
and $\eta_0=1$. Let $N_5$ be a positive integer satisfying
$\eta_{N_5}>\frac{2(\gamma+\theta)}{\gamma+\theta+1}-\eta$. It is
easy to see that $\eta<2$ and $\kappa_j<1$, $j=0,1\ldots,N_5$.
Using  similar arguments as that in Lemma \ref{ss-L3.14}, applying
the induction method, we can obtain the following lemma and omit
the proof.
\begin{lem}\label{ss-L3.15}
 Under the assumptions of Theorem \ref{ss-thm}, we have
  \begin{equation}
      \int^M_0\left\{u^2(x,t)+(M-x)^{\frac{\gamma-1}{\gamma}}(g-g_\infty)^2
      +(r-r_\infty)^2\right\}dx\leq
      C_{\eta,j}(1+t)^{-\eta_j},\label{ss-E5.1}
        \end{equation}
        \begin{equation}
    \int^t_0(1+s)^{\eta_j}\|u(\cdot,s)\|_{L^\infty}^2ds+    \int^t_0\int^M_0
    (1+s)^{\eta_j}\left(\rho^{\theta+1}u_x^2+\rho^{\theta-1}u^2
        \right)(x,s)dxds\leq C_{\eta,j},
        \end{equation}
        \begin{equation}
                \left\|\rho^\frac{\gamma+\theta}{2}(\cdot,t)-\rho_\infty^\frac{\gamma+\theta}{2}(\cdot)
                \right\|_{L^\infty}\leq
             C_{\eta,j}(1+t)^{-\frac{\eta_j}{4}},\label{ss-E5.2}
            \end{equation}
        \begin{equation}
   \int^{t}_0\int^M_0   (1+s)^{\kappa_j}
   \left[ \rho_\infty^{\gamma-1}
    (g-g_\infty)^2+(r-r_\infty)^2\right]dxds\leq
    C_{\eta,j},\label{ss-E5.3}
        \end{equation}
            \begin{equation}
             (1+t)^{\eta_j} \int^M_0\left(\rho^{\theta-1}
              u^2+\rho^{\theta+1}u_x^2
              \right)(x,t)dx+\int^t_0\int^M_0(1+s)^{\eta_j}u_t^2(x,s)dxds\leq
              C_{\eta,j},
            \end{equation}
                \begin{equation}
                  \|u(\cdot,t)\|_{L^\infty}\leq C_{\eta,j}
                  (1+t)^{-\frac{\eta_j}{2}},\label{ss-E5.10}
                \end{equation}
 for all $t\geq0$, $j=0,\ldots,N_5$, where $C_{\eta,j}$ is a
 positive constant depending on $\eta$ and $j$.

\end{lem}
\begin{rem} The main difficulty is to show (\ref{ss-E5.3}) with $j=k$, when
(\ref{ss-E5.1})-(\ref{ss-E5.2}) hold with $j=k$. From
(\ref{ss-E2.22-1})-(\ref{ss-E2.22-1-2}), we have
 \begin{eqnarray}
      &&\int^T_0\int^M_0(1+t)^{\kappa_{k}}
      \left[A(\rho_\infty^\gamma-\rho^\gamma)(\rho^{-1}-\rho_\infty^{-1})
      +G(M_0+x)(r^{2-2n}-r_\infty^{2-2n})(\frac{r^n}{n}-\frac{r^n_\infty}{n})\right]dxdt\nonumber\\
            &=&-\int^T_0\int^M_0(1+t)^{\kappa_{k}}\frac{u_t}{r^{n-1}}(\frac{r^n}{n}-\frac{r^n_\infty}{n})dxdt
            -\int^T_0\int^M_0(1+t)^{\kappa_{k}}\Delta fr^{1-n}\left(\frac{r^n}{n}-\frac{r^n_\infty}{n}
      \right)dxdt
            \nonumber\\
            &&+\int^T_0\int^M_0(1+t)^{\kappa_{k}}
            (2c_1+c_2)\rho^{1+\theta}\partial_x(r^{n-1}u)(\rho_\infty^{-1}-\rho^{-1})dxdt\nonumber\\
     &&+\int^T_0\int^M_0(1+t)^{\kappa_{k}}
            2c_1(n-1)\rho^\theta \left(\frac{u}{r}(\frac{r^n}{n}-\frac{r^n_\infty}{n})
                   \right)_xdxdt:=\sum^4_{i=1}Q_i,\ T>0\label{ss-E5.22-5}
    \end{eqnarray}
and
    \begin{equation}
    \textrm{L.H.S of } (\ref{ss-E5.22-5})\geq C_{12}\int^T_0\int^M_0(1+t)^{\kappa_{k}}
     \left[ \rho_\infty^{\gamma-1}(g-g_\infty)^2+
      \left(r-r_\infty\right)^2
      \right]dxdt.
    \end{equation}
Similar to (\ref{ss-E3.29})-(\ref{ss-E2.30-1}), applying the
estimates (\ref{ss-E2.9})-(\ref{ss-E2.13}), (\ref{ss-E3.92}),
 integration by parts, the Cauchy-Schwarz inequality and the
fact that
$\kappa_j=\frac{3\gamma+3\theta-1}{4(\gamma+\theta)}\eta_j-\frac{\epsilon_4}{4}\eta_j-\frac{1}{2}-\frac{\epsilon_4}{2}
<\eta_j$, we can estimate $Q_i$ as follows.
\begin{eqnarray}
  Q_1&\leq& -\left.\int^M_0(1+t)^{\kappa_{k}}\frac{u}{r^{n-1}}\left(\frac{r^n}{n}-\frac{r_\infty^n}{n}
      \right)dx\right|^T_0+
      C\int^T_0(1+t)^{\kappa_{k}}\|u\|_{L^\infty_x}^2dt\nonumber\\
            &&+C\int^T_0\int^M_0(1+t)^{\kappa_k-1}|u||r-r_\infty|dxdt\nonumber\\
            &\leq& C
\end{eqnarray}
    \begin{eqnarray}
      Q_2&\leq &\frac{C_{12}}{6}\int^T_0\int^M_0(1+t)^{\kappa_{k}}
      \left(r-r_\infty\right)^2
      dxdt+C\int^T_0f_1^2(1+t)^{\kappa_{k}}dt\nonumber\\
                &\leq &\frac{C_{12}}{6}\int^T_0\int^M_0(1+t)^{\kappa_{k}}
      \left(r-r_\infty\right)^2
      dxdt+C,
    \end{eqnarray}
    \begin{eqnarray}
      Q_3  &\leq&C \int^T_0\int^M_0(1+t)^{\eta_{k}}\rho^{1+\theta}(r^{n-1}u_x)_x^2dxdt\nonumber\\
      &&+C
      \int^T_0\int^M_0(1+t)^{2\kappa_{k}-\eta_{k}-\frac{\eta_{k}\nu_1}{2}}(M-x)^{\frac{\theta-1}{\gamma}
      -\frac{\nu_1(\gamma+\theta)}{\gamma}}dxdt\nonumber\\
        &\leq&C,
    \end{eqnarray}
where $\nu_1=\frac{\gamma+\theta-1}{\gamma+\theta}-\epsilon_4$,
    $$
    \|r(\cdot,t)-r_\infty(\cdot)\|_{L^\infty}\leq
    C(1+t)^{\frac{\nu_1\eta_{k}}{4}},
    $$
and
    \begin{eqnarray}
      Q_4&\leq&C \int^T_0\int^M_0(1+t)^{\eta_{k}}\left(\rho^{1+\theta}(r^{n-1}u_x)_x^2
      +\rho^{\theta-1}u^2
      \right)dxdt\nonumber\\
      &&+C \int^T_0\int^M_0(1+t)^{2\kappa_{k}-\eta_{k}-\frac{\eta_{k}\nu_1}{2}}(M-x)^{\frac{\theta-1}{\gamma}
      }dxdt\nonumber\\
        &\leq&C.\label{ss-E5.30-7}
    \end{eqnarray}
From  (\ref{ss-E5.22-5})-(\ref{ss-E5.30-7}),  we  finish the proof
of (\ref{ss-E5.3}) with $j=k$.
\end{rem}
From (\ref{ss-E5.3}) with $j=N_5$, using  similar arguments as
that in Lemma \ref{ss-L2.3}, we can obtain the following lemma and
omit the proof.
\begin{lem}
 Under the assumptions of Theorem \ref{ss-thm}, we have
    \begin{equation}
    \int^M_0(M-x)^{2-\frac{2\theta}{\gamma}}(\rho^\theta-\rho_\infty^\theta)_x^2dx\leq
    C_\eta(1+t)^{\eta-\frac{\gamma+\theta-1}{\gamma+\theta}},
    \end{equation}
        \begin{equation}
           \|\rho^{\gamma}(\cdot,t)-\rho_\infty^{\gamma}(\cdot)\|_{L^\infty}\leq
             C(1+t)^{\frac{\eta}{2}-\frac{3\gamma+3\theta-1}{4(\gamma+\theta+1)}},
        \  t\geq0.
        \end{equation}
\end{lem}
Thus, we finish the proof of Theorem \ref{ss-thm}.


\begin{thebibliography}{50}
\bibitem{Bresch2006}
    D. Bresch and B. Desjardins,
        {On the construction of approximate solutions for the 2D viscous shallow water model and for compressible Navier-Stokes
        models},
        J. Math. Pure Appl., 86(2006), pp. 362--368.
\bibitem{Chen2002}
      G.Q. Chen and M. Kratka,
    {Global solutions to the Navier-Stokes equations for compressible
    heat-conducting flow with symmetry and free boundary}, Comm.
    Partial Differential Equations, 27(5-6)(2002), pp. 907--943.
\bibitem{ChenP}
    P. Chen and T. Zhang,
    {A vacuum problem
    for multidimensional compressible Navier-Stokes equations with
    degenerate viscosity coefficients}, preprint,
    arXiv:math.AP/0701150.
\bibitem{Ducomet2005}
     B. Ducomet and A.A. Zlotnik,
    {Viscous compressible barotropic symmetric flows
    with free boundary under general mass force. I. uniform-in-time bounds and
    stabilization}, Math. Methods Appl. Sci., 28(7)(2005), pp.
    827--863.
\bibitem{Ducomet2005-2}
    B. Ducomet and  A.A. Zlotnik,
    {Lyapunov functional method for 1D radiative and reactive viscous gas
    dynamics},
    Arch. Rational Mech. Anal., 177(2005), pp. 185--229.
\bibitem{Grad}
    H. Grad,
    { Asymptotic theory of the Boltzmann equation II}, In:
     Rarefied Gas Dynamics, 1. J. Laurmann, ed., NewYork: Academic Press,
    1963, pp. 26--59.
\bibitem{Hartman}
    P. Hartman,
    { Ordinary differential equations},
    Wiley, New York, 1964.
\bibitem{hoff}
      D. Hoff and  D. Serre,
    {The failure of continuous dependence on initial data for the
    Navier-Stokes equations of compressible flow}, SIAM J. Appl.
    Math.,
    51(4)(1991), pp.   887--898.
\bibitem{Kobayashi}
    T. Kobayashi and Y. Shibata,
    {Decay estimates of solutions for
    the equations of motion of compressible viscous and
    heat-conductive gases in an exterior domain in $ R\sp 3$},
     Comm. Math. Phys., 200(3)(1999), pp. 621--659.

\bibitem{Lions}
     P.L.  Lions,
    { Mathematical topics in fluid mechanics}, Vol.
    \textbf{1-2},
    Oxford University Press: New York, 1996, 1998.
\bibitem{liu}
     T.P.  Liu, Z.P.  Xin and  T. Yang,
    {Vacuum states of compressible flow},
    Discrete Contin. Dynam. Systems, 4(1)(1998), pp.  1--32.
\bibitem{Matsumura}
    A. Matsumura and  S. Yanagi,
    {Uniform boundedness of the solutions for a one-dimensional isentropic
    model system of a compressible viscous gas},
    Comm. Math. Phys., 175(1996), pp. 259--274.
\bibitem{MatusuNecasova}
     \u{S}.  Matu\u{s}\.{u}-Ne\u{c}asov\'{a},   M. Okada and
      T. Makoni,
    {Free boundary problem for the equation of spherically symmetric
    motion of viscous gas (II)-(III)},
    Japan J. Indust. Appl. Math., 12(1995) pp. 195--203; 14(1997),
    pp.    199--213.
\bibitem{okada93}
   M. Okada  and  T. Makino,
    {Free boundary value problems for the equation of spherically symmetrical motion of viscous
    gas},
    Japan J. Appl. Math., 10(1993), pp. 219--235.
\bibitem{okada1}
  M.  Okada,
     \u{S}.  Matu\u{s}\.{u}-Ne\u{c}asov\'{a} and  T. Makino,
    {Free boundary problem for the equation of one-dimensional motion
    of compressible gas with density-dependent viscosity}, Ann. Univ.
    Ferrara Sez. VII (N.S.), 48(2002), pp. 1--20.
\bibitem{Straskraba}
   I. Stra\u{s}kraba  and  A.A.  Zlotnik,
    {Global behavior of 1d-viscous compressible barotropic fluid with a free boundary and large
    data},
    J. Math. Fluid Mech., 5(2003), pp. 119--143.
\bibitem{Vaigant}
    V.A.  Vaigant and A.V. Kazhikhov,
    {On existence of global solutions to the two-dimensional
    Navier-Stokes equations for a compressible viscosity fluid},
    Siberian Math. J., 36(1995), pp. 1108--1141.
\bibitem{Ukai}
   S.  Ukai, T. Yang and H.J. Zhao,
    {Convergence rate for the compressible Navier-Stokes equations with external
    force},
    J. Hyperbolic Differ. Equ., 3(3)(2006), pp. 561--574.
\bibitem{vong}
   S.W.  Vong,  T. Yang and  C.J. Zhu,
    {Compressible Navier-Stokes equations with degenerate viscosity
    coefficient and vacuum(II)}, J. Differential Equations,
    192(2)(2003), pp.      475--501.
\bibitem{xin}
  Z.P.  Xin,
    {Blow-up of smooth solution to the compressible Navier-Stokes
    equations with compact density}, Comm. Pure Appl. Math.,
    51(3)(1998), pp.      229--240.
\bibitem{yang3}
 T.  Yang   and C.J. Zhu,
    {Compressible Navier-Stokes equations with degenerate viscosity
    coefficient and vacuum},
     Comm. Math. Phys., 230(2)(2002), pp. 329--363.
\bibitem{fang06}
   T. Zhang and  D.Y. Fang,
    {Global behavior of  compressible
    Navier-Stokes equations with a degenerate viscosity
    coefficient},
    Arch. Rational Mech. Anal., 182(2)(2006), pp. 223--253.
\bibitem{fang06-2}
   T. Zhang and  D.Y. Fang,
    {Global behavior of  spherically symmetric
    Navier-Stokes equations with density-dependent viscosity}, preprint,
    arXiv:math.AP/0701216.
\bibitem{Zlotnik2005}
     A.A. Zlotnik   and B.  Ducomet,
      {The stabilization rate and stability of viscous compressible barotropic
      symmetric flows with a free boundary for a general mass
      force},
      Sb. Math., 196(11-12)(2005), pp. 1745--1799.
\end{thebibliography}
\end{document}